\documentclass{article}
\usepackage[utf8]{inputenc}
\pdfoutput=1

\usepackage{amsmath,amssymb,amsthm}
\usepackage{mathtools}

\usepackage{tikz}

\usepackage[backend=bibtex,url=false]{biblatex}
\usepackage{cleveref}

\usepackage{array} 
\newcolumntype{C}[1]{>{\centering\arraybackslash}p{#1}} 

\newtheorem{theorem}{Theorem}[section]
\newtheorem{definition}[theorem]{Definition}

\newtheorem{lemma}[theorem]{Lemma}
\newtheorem{corollary}[theorem]{Corollary}

\newtheorem{question}[theorem]{Question}
\newtheorem{remark}[theorem]{Remark}

\def\N{\mathbb{N}}
\def\Z{\mathbb{Z}}

\def\R{\mathbb{R}}

\newcommand{\dloc}{d_{\mathrm{loc}}}
\DeclareMathOperator{\cl}{cl}
\DeclareMathOperator{\Int}{int}

\usepackage{authblk}

\title{New examples of $M\setminus L$: intruder sets}

\author[1]{Harold Erazo}

\affil[1]{IMPA, Estrada Dona Castorina, 110. Rio de Janeiro, Brazil. \texttt{harold.eraz@gmail.com}}

\date{\today}

\begin{document}

\maketitle

\begin{abstract}
    We exhibit new examples of regions of $M\setminus L$ where $M$ and $L$ denote the Markov and Lagrange spectra, respectively. These regions have a different nature from all known regions studied so far: they contain \emph{intruder sets} associated with distinct combinatorics that trespass the region where self-replication holds. Our construction follows the usual self-replication method but replaces the standard local uniqueness condition with a more flexible and weaker property. These examples emerged from a large-scale computational search for regions of $M\setminus L$, which indicates that many such regions with intruder sets exist. We conclude with some open problems about these new regions.
\end{abstract}



\section{Introduction}

The \emph{Lagrange spectrum} consists of the best constants of Diophantine approximations of irrational numbers
\begin{equation*}
    L:=\left\{k(\alpha):=\limsup\limits_{\substack{p,q\to\infty\\p,q\in\mathbb{Z}}} \frac{1}{|q(q\alpha-p)|}<\infty: \alpha\notin\mathbb{Q}\right\}.
\end{equation*}
The \emph{Markov spectrum} is related to the approximation of binary quadratic forms and is defined to be
\begin{equation*}
    M:=\left\{\sup_{(p,q)\in\Z^{2}\setminus\{(0,0)\}}\frac{1}{|ap^{2}+bpq+cq^{2}|}\,:\,ax^{2}+bxy+cy^{2}\text{ real indefinite, }b^{2}-4ac = 1\right\}.
\end{equation*}
\[\]
It is known that $L\subset M\subset\R^{+}$ are closed subsets.

Both spectra have been studied since the works of Markov \cite{M:formes_quadratiques1} and \cite{M:formes_quadratiques2} from 1879 and 1880 showing that 
\[
    L \cap [0, 3) = M \cap [0, 3) = \left\{ \sqrt{5} < \sqrt{8} < \frac{\sqrt{221}}{5} < \dotsb \right\},
\]
that is, $L$ and $M$ coincide below $3$, and they consist of a discrete sequence of quadratic irrationals. In 1975, Freiman~\cite{Fr75} showed that $[c_{F},\infty)\subset L\subset M$ and $(\nu_{F},c_F)\cap M=\varnothing$ where $c_F = 4.52782956\ldots\in L$ and $\nu_F = 4.52782953\ldots\in L$. This ray $[c_{F},\infty)$ is known as \emph{Hall's ray} and $c_F$ is known as \emph{Freiman constant}. In particular the difference set $M\setminus L$ is a subset of $(3,\nu_F)$.

The first examples of elements in $M\setminus L$ were discovered in 1968 by Freiman \cite{F:non-coincidence}. Some years latter, in 1973 Freiman \cite{Fr73} found another point of $M\setminus L$. In 1977, Flahive \cite{Flahive77} proved that the second example of Freiman is accumulated by a discrete sequence of points in $M\setminus L$. The subject remained dormant until recently (around 2017) when C. Matheus and C. G. Moreira characterized completely in \cite{Freiman311} and \cite{Freiman329} the previous regions studied by Freiman, and as a by-product, showed that the difference set has positive Hausdorff dimension. Pushing further the method, the same authors \cite{MM:markov_lagrange} found the first example of a region of $M\setminus L$ above $\sqrt{12}$ (thus disproving a conjecture by Cusick). Further explicit examples, both above and below $\sqrt{12}$, have been studied in \cite{Newgaps}, \cite{RGH}, \cite{LGH}, \cite{MminusLisnotclosed}, \cite{MminusLnear3}. The complete list of all known explicit examples of regions of $M\setminus L$ is given in \Cref{tab:knownregions}.

The previous examples give evidence that the nature of the spectra between 3 and $c_F$ is rich. Indeed, in \cite{M:geometric_properties_Markov_Lagrange}, it was shown that the Hausdorff dimensions of both spectra are equal if we truncate them, more precisely, one has for any $t\in\R$ that
\begin{equation*}
    d(t):=\dim_H\left(L\cap(-\infty,t]\right)=\dim_H\left(M\cap(-\infty,t]\right),
\end{equation*}
and that the function $t\mapsto d(t)$ is continuous and such that $d(3+\varepsilon)>0$ for all $\varepsilon>0$. Moreover, the function
\begin{equation*}
    D(t):=\dim_H\left(k^{-1}(-\infty,t)\right)=\dim_H\left(k^{-1}(-\infty,t]\right)
\end{equation*}
is also continuous and satisfies $d(t)=\min\{1,2\cdot D(t)\}$. 

On the other hand, the difference set $M\setminus L$ is a non closed set \cite{MminusLisnotclosed} that is still not well understood. The best bounds for its Hausdorff dimension are (see \cite{MM:markov_lagrange,LGH,MMPV}), 
\begin{equation*}
    0.594561<\dim_H\left(M\setminus L\right)<0.796445.
\end{equation*}
In particular, the interiors of the two spectra coincide, $\Int(L)=\Int(M)$, so the difference set $M\setminus L$ is contained in the complement of this common interior. Since it is unknown whether $L\cap(-\infty,\nu_F]$ has nonempty interior, locating elements of the difference set is especially challenging. The earliest appearance of this question can be tracked back to Berstein’s 1973 conjecture \cite[Page 77]{BerConj}, who conjectured that $[4.1,4.52]\subset L$.

In \cite[Page 149]{M:geometric_properties_Markov_Lagrange}, it was acknowledged that describing the geometric structure of the difference set $M\setminus L$ is a central question to understand the spectra. In this paper, we present a detailed study of two new examples of $M\setminus L$ whose geometry is significantly more intricate than that of all previously known regions. These examples emerged from an extensive computational search, the details of which will be provided in a forthcoming publication, but since the nature of the intruders sets is interesting by itself, we decided to create this separate paper. Additionally, we compile a list of numerous further examples that display the same complex structure, adding to the growing evidence that the geometry of $M\setminus L$ is considerably richer than previously understood. 

The remainder of the paper is organized as follows. In Section 2, we review the previously known regions of the difference set $M\setminus L$ and contrast them with the new regions that contain intruder sets. Section 3 introduces additional definitions and notation used throughout the paper. In Section 4, we provide a complete analysis of the first example in the difference set that includes an intruder set. Section 5 offers a partial analysis of a second such example. In Section 6, we establish a bound on a region of the spectra with small Hausdorff dimension. Finally, in \Cref{sec:further_examples}, we present many additional orbits that exhibit intruder sets.

\noindent\textbf{Acknowledgements:} I would like to express my gratitude to Carlos Gustavo Moreira for many helpful discussions about this paper. The author is partially supported by CAPES and FAPERJ.

\noindent\textbf{Data Availability:} The datasets containing the orbits reported in \Cref{tab:alphabet_12,tab:alphabet_123,tab:regions_small_dimension} and the algorithms used in this paper are available at \cite{repo}.


\section{Description of the new examples and open problems}

\subsection{Dynamical definition of the Lagrange and Markov spectra}

To understand the nature of the new examples we will present, we need to introduce the dynamical characterization of the spectra through symbolic dynamics first. 

Let $\Sigma := (\N^*)^{\Z}$ be the space of bi-infinite sequences and $\sigma:\Sigma\to\Sigma$ the left-shift map $\sigma((a_n)_{n\in\Z}):=(a_{n+1})_{n\in\Z}$. When referring to orbits we always refer to orbits of the dynamical system $(\Sigma,\sigma)$. 

For $k\in\Z$ define the function  $\lambda_k : \Sigma \rightarrow \R$ by
\begin{equation*}
    \lambda_k((a_n)_{n\in \Z}) = [a_k;a_{k+1},\dots]+[0;a_{k-1},a_{k-2},\dots]
\end{equation*}
where $\left[b_0,b_1,b_2,\dots\right] = b_0+\frac{1}{b_1+\frac{1}{b_2+\dots}}$ is the usual continued fraction associated with the sequence of positive integers $(b_n)_{n\geq 0}$. The \emph{Lagrange value} of a sequence $\underline{a}=(a_n)_{n\in\Z}\in\Sigma$ is given by
\begin{equation*}
    \ell(\underline{a}):=\limsup_{k\to\infty}\lambda_k(\underline{a})
\end{equation*}
and the \emph{Markov value} of $\underline{a}$ is
\begin{equation*}
    m(\underline{a}):=\sup_{k\in\Z}\lambda_k(\underline{a}).
\end{equation*}
In 1921, Perron proved that the Lagrange and Markov spectra can be characterized as
\begin{equation*}
    L:=\{\ell(\underline{a})<\infty:\underline{a}\in\Sigma\} \quad \text{and} \quad M:=\{m(\underline{a})<\infty:\underline{a}\in\Sigma\}.
\end{equation*}

It is not difficult to show that if $t\in M$, then there is a bi-infinite sequence $\underline{a}\in\Sigma$ such that $t=m(\underline{a})=\lambda_0(\underline{a})$. 

Given an infinite sequence $(a_{i})_{i\in\Z}\in\Sigma$, we write it as $\ldots a_{-2}a_{-1}a_{0}^{*}a_{1}a_{2}\ldots$ where the asterisk denotes the 0th position. We will also use an overline to denote periodicity, for example $\overline{1^{*}23} = \ldots1231^{*}23123\ldots$. Such a bi-infinite sequence is a \emph{periodic orbit}. Markov values of periodic orbits naturally belong to the Lagrange spectrum $L$. When writing an equation of the form $t=m(\dots a_{-1}a_0^*a_1\dots)$, we always refer that the Markov value is being attained at the marked position, that is $t=m(\dots a_{-1}a_0^*a_1\dots)=\lambda_0(\dots a_{-1}a_0^{*}a_1\dots)$.

Given an infinite sequence $\underline{a}=(a_i)_{i\in\Z}\in\Sigma$, the sequence $\underline{a}^T:=(a_{-i})_{i\in\Z}$ is called the \emph{transpose} of $\underline{a}$. Analogously, for finite words $a=a_1\dots a_n\in(\N^*)^n$ the transpose is $a^T:=a_n\dots a_1$. Given a marked finite word $w^*=a_{-p}\dots a_0^*\dots a_q\in\{1,\dots,A\}^{[-p,q]\cap\Z}$ and some $\underline{b}=(b_n)_{n\in\Z}\in\{1,\dots,A\}^\Z$, we write $\underline{b}=a_{-p}\dots a_0^*\dots a_q$ to denote that $(b_{-p},\dots,b_q)=(a_{-p},\dots,a_q)$. Given some $t\in\R$, an inequality of the form $\lambda_0(w^*)>t$ means that
\begin{equation*}
    \inf\left\{\lambda_0(\underline{b}):\underline{b}\in\{1,\dots,A\}^Z,\underline{b}=w^*\right\}=\min_{u,v\in\left\{\overline{1A},\overline{A1}\right\}}\lambda_0(u\,w^*\,v)>t. 
\end{equation*}
In other terms, $\lambda_0(w)>t$ means that any extension of $w=a_{-p}\dots a_0^*\dots a_q$ in the alphabet $\{1,\dots,A\}$ gives Markov values greater than $t$, where the alphabet will be always clear from context. Analogous definitions of all the above notations are used for one-sided words.

\subsection{Geometry of classical regions of $\boldsymbol{M\setminus L}$}

The usual construction of regions in $M\setminus L$ relies on the concepts of \emph{non semisymmetric words} (\Cref{def:nonsemisymmetric}), \emph{local uniqueness} (\Cref{def:local_uniqueness}) and \emph{self-replication} (\Cref{def:self-replication}). This method, developed in \cite{MM:markov_lagrange}, remains the only known technique for generating elements in the difference set. In fact, for certain regions of the spectra, it is known that it is the only possible construction (see \cite[Proposition 2.16]{LGH}).

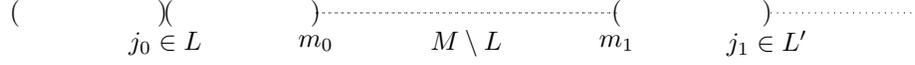
\begin{figure}[ht]
    \centering
    \begin{tikzpicture}[scale=4]
        \node[] at (-0.5,0) {(}; 
        \node[] at (-0.01,0) {)}; 
        \node[] at (0.01,0) {(}; 
        \draw[dotted] (1.5,0) -- (0.5,0) node {)};
        \draw[dotted] (0.5,0) -- (1.5,0) node {(};
        \draw[-,dotted] (2.5,0) -- (2,0) node {)};
        \node[] at  (0,-0.1) {$j_0\in L$};
        \node[] at  (2,-0.1) {$j_1\in L^\prime$};
        \node[] at  (0.5,-0.1) {$m_0$};
        \node[] at  (1.5,-0.1) {$m_1$};
        \node[] at  (1,-0.1) {$M\setminus L$};
    \end{tikzpicture}
    \caption{Structure of known regions of $M\setminus L$. The interval $(j_0,j_1)$ is a gap of $L$ containing points of $M$. The set $L^\prime$ corresponds to the accumulation points of $L$.}\label{fig:MminusLclassic}
\end{figure}

Given a set $X\subset \R$ we denote by $X^\prime$ the set of accumulation points of $X$. All known regions of $M\setminus L$ have the structure depicted in Figure~\ref{fig:MminusLclassic} where the spaces represent gaps in $M$, $j_1\in L^\prime$, $j_0$ is isolated in $M$, and $(M\setminus L)\cap(j_0,j_1)$ is given as a disjoint union of the form
\begin{equation}\label{eq:mldecomposition}
    (M\setminus L)\cap(j_0,j_1)=C\sqcup D
\end{equation}
where $C$ is a Cantor set and $D$ is an infinite set of discrete points in $M$. Moreover for such regions it holds that $D^\prime=C$.

Such regions are usually associated with a odd non semisymmetric word $w$. The left boundary $j_0:=m(\overline{w})=\ell(\overline{w})$ corresponds to the Markov value of the periodic orbit with period $w$. For all known regions this word is unique\footnote{For all known regions of $M\setminus L$ this is the case. In general, there can be a finite set of non semisymmetric distinct orbits $\overline{w_1},\overline{w_1^T},\dots,\overline{w_r},\overline{w_r^T}$ such that $j_0 = m(\overline{w_i})=m(\overline{w_i^T})$ for all $i=1,\dots,r$.} up to cyclic permutations and transpositions, because for these regions the unique orbit with Markov value $j_0$ is precisely $\overline{w}$. The right boundary $j_1=j_1(w)$ corresponds to the \emph{self-replication height} (see \Cref{def:self-replication-height}), which is the minimal Markov value such that self-replication holds. All sequences $m(\underline{a})=\lambda_0(\underline{a})\in(M\setminus L)\cap(j_0,j_1)$ connect to $w$ to some side up to transposition, that is, either $\underline{a}$ or $\underline{a}^T$ is eventually periodic to the left or to the right with period $w$. In fact, writing $w=w_1\dots w_n$, for all such regions it holds that for some $m\in\N^*$, the set
\begin{equation}\label{eq:setP}
    P = \{m(B)=\lambda_0(B)\in(j_0,j_1): B=\overline{w}w^*w_1\dots w_m\ldots\in\{1,\dots,A\}^\Z\}
\end{equation}
is equal to $(M\setminus L)\cap(j_0,j_1)$ (with the convention $w_{k}:=w_{k-n}$ when $k>n$). The fact that $j_1\in L^\prime$ for all known regions is a consequence of \Cref{lem:flahive}.

For all these regions, one can characterize the sets appearing in \eqref{eq:mldecomposition} as follows:
\begin{equation}\label{eq:setC}
    C = \left\{m(\underline{a})=\lambda_0(\underline{a})\in(j_0,j_1):\underline{a}\in(\N^*)^\Z, \min\{\ell(\underline{a}),\ell(\underline{a}^T)\}<j_0=\max\{\ell(\underline{a}),\ell(\underline{a}^T)\}\right\},
\end{equation}
\begin{equation}\label{eq:setD}
    D = \left\{m(\underline{a})=\lambda_0(\underline{a})\in(j_0,j_1):\underline{a}\in(\N^*)^\Z, \ell(\underline{a})=\ell(\underline{a}^T)=j_0\right\}.
\end{equation}
Moreover, the sequences $m(\underline{a}) = \lambda_0(\underline{a})$ taking values in $C$ can be equivalently characterized as those that connect (up to transposition) with $w$ on exactly one side, while the sequences taking values in $D$ correspond to those that connect with $w$ on one side and with $w^T$ on the other. 

Observe that according to \Cref{fig:MminusLclassic}, for all known regions of $M\setminus L$, it holds that $j_0=\ell(\overline{w})$ is isolated in $M$ (so also in $L$) and that if $(j_0,r_0)$ is the maximal gap of $L$ to the right of $j_0$, i.e., $L\cap(j_0,r_0)=\varnothing$ with $j_0,r_0\in L$, it holds that $(M\setminus L)\cap(j_0,r_0)\neq\varnothing$.

Finally, for all known regions of $M\setminus L$, it holds that the local Hausdorff dimension of $M$ 
\begin{equation*}
    \dloc^M(t):=\lim_{\varepsilon\to}\dim_H\left(M\cap(t-\varepsilon,t+\varepsilon)\right)
\end{equation*}
is constant in the interval $(j_0,j_1)$ (see \cite[Theorem 1.2]{LGH}).
{
\renewcommand{\arraystretch}{1.4}
\begin{table}[h!]
    \begin{tabular}{C{3.1cm}C{1cm}C{2.5cm}C{3.5cm}C{2.7cm}}
        \hline
        Non semisymmetric word & Length & Markov value of periodic orbit & Discovery & Complete description \\ 
        \hline
        122122221 & 9 & $3.1181201\dots$
        & Freiman, 1968 \cite{F:non-coincidence} & \cite{Freiman311} \\
        2112122 & 7 & $3.2930442\dots$ & Freiman, 1973 \cite{Fr73} & \cite{Freiman329} \\
        1233222 & 7 & $3.7096998\dots$ & C.Matheus, C.G.Moreira, 2020 \cite{MM:markov_lagrange} & \cite[Section 2.5]{CDMLS} \\
        $2^{2k}112^{2k+1}112^{2k+2}11$ $k=1,2,3,4$ & $6k+9$ & Near 3 & D.Lima, C.Matheus, C.G.Moreira, S.Vieira, 2021 \cite{MminusLnear3} & \\
        $2^{2k-1}12^{2k}12^{2k+1}1$ $k\geq 1$ & $6k+3$ & Near $1+\frac{3}{\sqrt{2}}=3.1213\dots$ & D.Lima, C.Matheus, C.G.Moreira, S.Vieira, 2022 \cite{MminusLisnotclosed} &  \\
        21133311121 & 11 & $3.9387762\dots$ & L.Jeffreys, C.Matheus, C.G.Moreira, 2024 \cite{Newgaps} & \cite[Section 4.5]{LGH} \\
        212332111 & 9 & $3.6766994\dots$ & C.Rieutord, C.G.Moreira, H.Erazo, 2024 \cite{RGH} & \cite{RGH} \\
        123332112 & 9 & $3.6911783\dots$ & C.Rieutord, C.G.Moreira, H.Erazo, 2024 \cite{RGH} & \cite{RGH} \\
        12111233311133232 & 17 & $3.9420011\dots$ & H.Erazo, L.Jeffreys, C.G.Moreira, 2024 \cite{LGH} & \cite[Theorem 4.43]{LGH} \\ \hline
    \end{tabular}
    \caption{Complete list of all known regions of $M\setminus L$ before our algorithmic search. The full characterization of some regions was done in posterior work, which is pointed out in the last column. The examples near $3$ and $1+\frac{3}{\sqrt{2}}$ were not completely characterized, but it is easy to do it and show that they also have the structure of \Cref{fig:MminusLclassic}.}
    \label{tab:knownregions}
\end{table}
}

\subsection{Regions of $M\setminus L$ with intruder sets}

In this paper we study new examples of $M\setminus L$ that differ from all previously known regions because they have some ``intruder sets" (see \Cref{fig:MminusLintruder1}). For these new regions, some of the properties described in the previous subsection are no longer true. For example, for these new regions, the orbit $j_0=m(\overline{w})$ is also isolated in $M$, but the maximal gap in $L$ to the right does not intersect $M$, so is also a maximal gap of $M$. Moreover, for these new regions, the simple characterization of the sets $C$ and $D$ given in \eqref{eq:setC} and \eqref{eq:setD} does not hold.

Recall from \eqref{eq:setP}, that for the regions of \Cref{tab:knownregions}, all sequences $m(\underline{a})=\lambda_0(\underline{a})\in(M\setminus L)\cap(j_0,j_1)$ share the fact that central portions of them are equal to the orbit $\overline{w}$. In fact, they connect (up to transposition) to the left or to the right with $\overline{w}$. For regions with intruders sets, there is a sum of Cantor sets $K_1+K_2$ that is a subset of $M\cap(j_0,j_1)$ but that is associated with a \emph{completely different} combinatorics from that of $\overline{w}$. In other words, the term intruder is motivated because of the fact that sequences that belong to intruder sets $m(\underline{a})=\lambda_0(\underline{a})\in K_1+K_2$ will have a different expression from those of $P$ and also contain points of $L^\prime$ (see \Cref{fig:MminusLintruder1} for an example). We call these sums $K_1+K_2$ as \emph{intruder sets}.

The two regions of $M\setminus L$ studied in this paper exhibit intruder sets. Moreover, these regions all lie in \emph{good intervals} (see \cite[Theorem 3.1]{LGH}). Hence, we have a good understanding of the local Hausdorff dimension of both spectra in the region $(j_0,j_1)$. In particular, from \cite[Theorem 1.2]{LGH} and \cite[Proposition 1.6]{LGH}, the local dimension $\dloc^M$ of $M$ is not constant in the interval $(j_0,j_1)$, since $D(t)$ is forced to increase when the parameter $t$ passes through points of $L^\prime$ in the intruder sets.

The regions studied in this paper are associated with the words $w=11122322221$ (\Cref{sec:onegap}) and $w=1112122$ (\Cref{sec:length7}).

\subsubsection{An example with alphabet $\{1,2,3\}$}\label{subsec:first_example}

Let $w=11122322221$ and $w^*=111223^*22221$. We mark precisely that position because 
\begin{equation*}
    j_0:=m(\overline{w})=\lambda_0(\overline{111223^*22221})=3.8340731702358\dots.
\end{equation*}
The Markov value
\begin{align*}
    j_1:=m(\overline{21}1232321212112221www^*w11121\overline{12})&=3.83407317023586722\ldots \\
    &\approx j_0+2.3728\cdot10^{-14}.
\end{align*}
is the self-replication height of $w$. Any bi-infinite sequence with Markov value below $j_1$ that is sufficiently similar to the periodic orbit $\overline{w}$ must self-replicate to one side (in this case is the right), i.e., must connect with $w$ (see \Cref{def:self-replication}). In particular the set
\begin{equation*}
    P:=\{m(B)=\lambda_0(B)\in(j_0,j_1): B=\dots2221ww^*\overline{w}\in\{1,2,3\}^\Z\},
\end{equation*}
corresponds to the natural candidates for points in $(M\setminus L)\cap(j_0,j_1)$. Unlike the simple regions described in the previous section, for regions that exhibit intruder sets, one only expects that $P\subset(M\setminus L)\cap(j_0,j_1)$ and the containment should be \emph{strict}. This is because of the existence of points in $L^\prime\cap(j_0,j_1)$, that can potentially produce more points of $M\setminus L$ associated with distinct non semi-symmetric orbits (with completely different combinatorics).

\begin{theorem}
    We have that $P\subset (M\setminus L)$ and $L^\prime\cap(j_0,j_1)\neq\varnothing$.
\end{theorem}

\begin{figure}[ht]
    \centering
    \begin{tikzpicture}[scale=4]
        \node[] at (-0.5,0) {(}; 
        \node[] at (-0.01,0) {)}; 
        \node[] at (0.01,0) {(}; 
        \draw[dotted] (0.3,0) -- (0.7,0) node {(};
        \draw[dotted] (0.7,0) -- (0.3,0) node {)};
        \draw[dotted] (1.5,0) -- (1,0) node {)};
        \draw[dotted] (1,0) -- (1.5,0) node {(};
        \draw[-,dotted] (2.5,0) -- (2,0) node {)};
        \node[] at  (0,-0.1) {$j_0\in L$};
        \node[] at  (2,-0.1) {$j_1\in L^\prime$};
        \node[] at  (0.3,0.1) {$i_0\in L^\prime$};
        \node[] at  (0.7,0.1) {$i_1\in L^\prime$};
        \node[] at  (1,-0.1) {$m_0$};
        \node[] at  (1.5,-0.1) {$m_1$};
        \node[] at  (1.25,0.1) {$M\setminus L$};
    \end{tikzpicture}
    \caption{Structure of $M$ near $j_0=m\left(\overline{111223^*22221}\right)=3.834\ldots$. The spaces represent maximal gaps of $M$. The interval $(i_1,j_1)$ is a gap of $L$ that contains points of $M$. There is a sum of Cantor sets $K_1+K_2\subset L\cap[i_0,i_1]$ that is ``trespassing" the region where self-replication holds.}
    \label{fig:MminusLintruder1}
\end{figure}
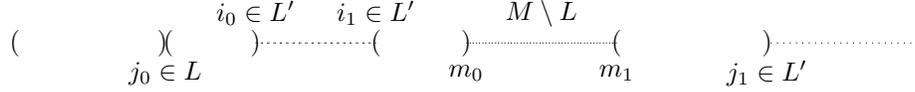

Let $m_0:=\min P$ and $m_1:=\max P$. In \Cref{thm1:characterization}, we show that in fact $(M\setminus L)\cap[m_0,m_1]=P$ has the usual decomposition: some Cantor sets and infinite isolated points that can be explicitly described. Moreover $L\cap(j_0,j_1)= L\cap[i_0,i_1]$ with $i_0,i_1\in L^\prime$ and $i_1<m_0$. The relation between these constants is shown in \Cref{fig:MminusLintruder1} and their explicit values are given in \Cref{sec:onegap}.

Observe that we do not claim that $(M\setminus L)\cap(i_0,i_1)=\varnothing$. However, if such points exists, they are associated with sequences that connect to non semisymmetric periodic orbits different from $\overline{11122322221}$. This is consequence of the results in \cite[Subsection 2.4]{LGH} and the fact that $(M\setminus L)\cap(i_0,i_1)$ is contained in a good interval.

\subsubsection{An example with alphabet $\{1,2\}$}

All orbits of period at most 5 in the alphabet $\{1,2\}$ are semisymmetric. The only non semisymmetric orbits of period 7 in the alphabet $\{1,2\}$ are precisely 2112122 and 1112122. Interestingly, these two orbits were already considered by Berstein in 1973 (see equation (2.11.5) in \cite[Page 45]{Berstein}). In that same page, Berstein proved that these two orbits satisfy the $\varepsilon$-property (see \cite[Page 42]{Cusick-Flahive}), which, in certain regions, is a necessary condition to construct elements of $M\setminus L$. The first word 2112122 was used by Freiman in 1973 to prove that $M\neq L$, thus producing the second example of elements in the difference set. In this paper we will show that the other orbit 1112122 also yields points in $M\setminus L$, however this requires much more effort than in the classical regions, since the existence of intruders sets yields points of $L^\prime$ very close to the points in $M\setminus L$. This could explain the fact why neither Freiman nor Berstein managed to use the other word 1112122 to produce more examples.

Consider the non semisymmetric word $w=1112122$.  and denote $w^*=1112^*122$. Define the constants 
\begin{equation*}       
    j_0:=m\left(\overline{w^*}\right)=m\left(\overline{1112^*122}\right)=3.334156277015\ldots.
\end{equation*}
and
\begin{align*}
    j_1:=m\left(\overline{212111112122}w^*ww121112\overline{221211111212}\right)&=3.3341574813\ldots \\
    &\approx j_0+1.204\cdot10^{-6}.
\end{align*}

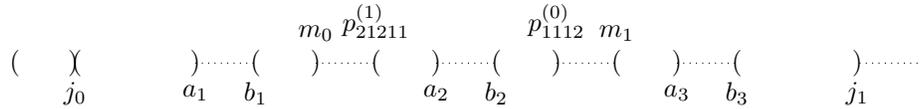
\begin{figure}[ht]
    \centering
    \begin{tikzpicture}[scale=0.8]
        \coordinate (j0) at (0,0);
        \coordinate (a1) at (2,0);
        \coordinate (b1) at (3,0);
        \coordinate (m0) at (4,0);
        \coordinate (m00) at (5,0);
        \coordinate (a2) at (6,0);
        \coordinate (b2) at (7,0);
        \coordinate (m11) at (8,0);
        \coordinate (m1) at (9,0);
        \coordinate (a3) at (10,0);
        \coordinate (b3) at (11,0);
        \coordinate (j1) at (13,0);

        \node at (-1,0) {(};
        \node at (-0.03,0) {)};
        \node at (0.03,0) {(};
        \node at (a1) {)};
        \node at (b1) {(};
        \node at (a2) {)};
        \node at (b2) {(};
        \node at (a3) {)};
        \node at (b3) {(};
        \node at (m0) {)};
        \node at (m00) {(};
        \node at (m11) {)};
        \node at (m1) {(};
        \draw[dotted] (a1) -- (b1) node {};
        \draw[dotted] (a2) -- (b2) node {};
        \draw[dotted] (a3) -- (b3) node {};
        \draw[dotted] (m0) -- (m00) node {};
        \draw[dotted] (m11) -- (m1) node {};
        \draw[dotted] (14,0) -- (j1) node {)};

        \node at (j0 |- 0,-0.5)  {$j_0$};
        \node at (a1 |- 0,-0.5)  {$a_1$};
        \node at (b1 |- 0,-0.5)  {$b_1$};
        \node at (m0 |- 0,0.5)  {$m_0$};
        \node at (m00 |- 0,0.7)  {$p_{21211}^{(1)}$};
        \node at (a2 |- 0,-0.5)  {$a_2$};
        \node at (b2 |- 0,-0.5)  {$b_2$};
        \node at (m11 |- 0,0.7)  {$p_{1112}^{(0)}$};
        \node at (m1 |- 0,0.5)  {$m_1$};
        \node at (a3 |- 0,-0.5)  {$a_3$};
        \node at (b3 |- 0,-0.5)  {$b_3$};
        \node at (j1 |- 0,-0.5)  {$j_1$};
    \end{tikzpicture}
    \caption{Partial structure of $M$ near $j_0=m\left(\overline{1112^*122}\right)=3.33415\ldots$ The elements $a_i,b_i\in L^\prime$, $i=1,2,3$ correspond to intruder combinatorics. The elements $m_0:=\min P$ and $m_1:=\sup P$ belong to $M\setminus L$ and one also has that $M\cap\left([m_0,p_{21211}^{(1)}]\cup[p_{1112}^{(0)},m_1]\right)\subset M\setminus L$, which are points from $P$. However, we do not know if $P\cap(a_2,b_2)\subset M\setminus L$.}
    \label{fig:complicated_intruder}
\end{figure}

The set
\begin{equation*}
    P:=\{m(B)=\lambda_0(B)\in(j_0,j_1): B=\overline{w}w^*w12\ldots\in\{1,2\}^\Z\},
\end{equation*}
corresponds to the natural candidates for points in $(M\setminus L)\cap(j_0,j_1)$. In \Cref{thm3:characterization1} and \Cref{thm3:characterization2} we will characterize $P\cap\left((b_1,a_2)\cup(b_2,a_3)\right)$ and we will show that

\begin{theorem}
    $P\cap\left((b_1,a_2)\cup(b_2,a_3)\right)\subset M\setminus L$.
\end{theorem}

Numerically, it is possible to verify that there are many more points in $P\cap(a_2,b_2)$ that also belong to $M\setminus L$. However, we were unable to decide if all of them belong to the difference set because of the existence of many intruder sets (which contain Cantor sets in $L^\prime$), which seem to ``intersect" $P$ in very small scales. In particular, we do not know if it is possible to characterize the set $P$ explicitly as was done in the example of \Cref{subsec:first_example}. 

\begin{question}\label{ques:main}
    Is it true that $P\subset (M\setminus L)\cap(j_0,j_1)$?
\end{question}

\subsection{Further examples of regions with intruders}

The examples discussed in the previous subsection are the result of a large-scale computer search for regions in $M\setminus L$. Many additional examples of regions in the difference set also exhibit intruder sets (see \Cref{sec:further_examples}). The orbits we chose to study explicitly in this paper are not particularly distinguished; they were selected simply because they have the smallest possible periods. Most of them display the same intricate behavior as the earlier examples: either the intruder sets are completely separated from $P$ as in \Cref{fig:MminusLintruder1} or they are so close to $P$ as in \Cref{fig:complicated_intruder}, that we are likewise unable to resolve the corresponding instance of \Cref{ques:main}.

Let us make precise the range of our computer explorations. An important parameter of transition of geometric complexity in the spectra is \[t_1 := \min\{t\in\R : d(t) = 1\}.\] In \cite{MMPV}, the authors confirmed previous heuristic estimates of Bumby and proved that
\begin{equation*}
    3.334384<t_1<3.334385.
\end{equation*}

Indeed, it was conjectured in
\cite[Page 149]{M:geometric_properties_Markov_Lagrange} that $t_1=\inf\Int L$, so it is conjectured that $t_1$ is being accumulated to the right by intervals of $L$, while $L$ is a fractal set to the left of $t_1$ (since it has dimension strictly smaller than 1). 

In particular, since the difference set lies outside the interior of the spectra, we expect that many orbits with Markov value greater than $t_1$ do not satisfy local uniqueness and therefore do not contribute to regions of $M\setminus L$. For this reason, our computational search has focused on orbits with Markov value less than $t_1$.

Let us call $s_0:=3.334384$ the lower bound for $t_1$. Recall that $\underline{a}\in(\N^*)^\Z$ has Markov value $m(\underline{a})\leq\sqrt{12}=3.46\ldots$ if and only if $\underline{a}\in\{1,2\}^\Z$. In particular, orbits below $t_1$ only contain $\{1,2\}$. We have tested all non semisymmetric orbits with odd minimal period less or equal than 23 in the alphabet $\{1,2\}$, and with Markov value at most $s_0$. Although most of these orbits generate regions of $M\setminus L$ in the classical manner (see \Cref{fig:MminusLclassic}), a significant number appear likely to produce regions containing intruder sets  (see \Cref{tab:alphabet_12}). We have verified computationally that some of these orbits indeed give rise to intruder sets, and we list several such examples in \Cref{subsec:intruder_alphabet_12}.

\begin{table}[h!]
\centering
    \begin{tabular}{C{1cm}C{4cm}C{3.5cm}}
        \hline
        Minimal period & Number of non semisymmetric orbits in $\{1,2\}$ before $s_0$ &  Number of orbits before $s_0$ that possibly exhibit intruder sets \\ 
        \hline
        7 & 2 & 1  \\
        9 & 9 & 1  \\
        11 & 36 & 4  \\ 
        13 & 140 & 21 \\
        15 & 502 & 68 \\
        17 & 1742 & 271 \\
        19 & 5964 & 978 \\
        21 & 20311 & 3473 \\
        23 & 69157 & 13672 \\
        \hline
    \end{tabular}
    \caption{Number of non-semisymmetric orbits (over the alphabet $\{1,2\}$) whose minimal period is odd and whose Markov value satisfies $\leq s_0=3.334384$, together with the number of these orbits that \emph{may} contain regions of $M\setminus L$ with intruder sets. We confirmed that the orbits listed in the third column with periods 7, 9, 11, 13 exhibit such regions (see \Cref{sec:further_examples}) and we can confirm these are the only ones for those periods 7-13.}
    \label{tab:alphabet_12}
\end{table}

We have also carried out a similar analysis in the alphabet $\{1,2,3\}$. A natural upper bound for our computational search is  $\nu_F = 4.52782953\ldots$, as points in the difference set lie entirely outside Hall’s ray  $[c_F,\infty)\subset L$, and are also clearly disjoint from Freiman’s gap $(\nu_F,c_F)$. Note also that any periodic orbit containing the digit 3 necessarily has Markov value at least $\sqrt{13}>t_1$. In \cite{Newgaps}, an analogous computer investigation was done to find the largest known elements in the difference set.

We have tested orbits in the alphabet $\{1,2,3\}$ with Markov value in the interval $[\sqrt{13},\nu_F]$ and with minimal period at most 15. For such orbits, the likelihood of generating regions of $M\setminus L$ is lower, since they may lie within connected components of the interior of the spectra. Nevertheless, we have discovered new examples of regions of $M\setminus L$ above $\sqrt{13}$ that also contain intruder sets. These orbits are listed in \Cref{subsec:intruder_alphabet_123}.

\begin{table}[h!]
\centering
    \begin{tabular}{C{1cm}C{5cm}C{5cm}}
        \hline
        Minimal period & Number of non semisymmetric orbits in $\{1,2,3\}$ in $[\sqrt{13},\nu_F]$ & Number of those orbits with confirmed intruder sets \\
        \hline
        5 & 12 & 0 \\
        7 & 115 & 0 \\
        9 & 932 & 0 \\
        11 & 7323 & 4 \\
        13 & 57447 & 10 \\
        15 & 448282 & 9 \\
        \hline
    \end{tabular}
    \caption{Number of non semisymmetric orbits with odd length minimal period over the alphabet $\{1,2,3\}$ whose Markov value lies between $\sqrt{13}$ and $\nu_{F} = 4.52782953\ldots$, and the number of such orbits for which we can confirm the existence of points of $M\setminus L$ with intruder sets. These orbits are explicitly listed in \Cref{subsec:intruder_alphabet_123}. We do not claim this list to be exhaustive.}
    \label{tab:alphabet_123}
\end{table}

\subsection{Regions of $M\setminus L$ with small dimension}

It is of particular interest to identify regions of $M\setminus L$ that contain intruder sets while having small Hausdorff dimension. We postpone a heuristic discussion of why such examples are relevant to the context of \Cref{ques:main} until the end of this section. First, we clarify what we mean by ``small dimension". Denote by $t_{1/3}$ the positive real number \[t_{1/3}:=\min\{t\in\R:D(t)=1/3\}.\] 
We are interested in searching regions of $M\setminus L$ in $M\cap(3,t_{1/3})$ that exhibit intruder sets. First, we need an explicit estimate of where this regime ends:

\begin{theorem}
\begin{equation*}
    \nu_1=3.037002\dots<t_{1/3}<\mu_2=3.037311\dots,
\end{equation*}
where $\nu_1$ and $\mu_2$ are defined in \eqref{eq:nu_1} and \eqref{eq:mu_2}, respectively.
\end{theorem}

However, after an broad computer search (see \Cref{tab:regions_small_dimension}), we were unable to find any example with intruder sets and with Markov values below $t_{1/3}$. This search was realized using an algorithm that we will explain in detail in a forthcoming work, which is partly based on the algorithm presented in \cite[Appendix B]{Newgaps}. From the exploration of all non semisymmetric orbits of odd minimal period at most 23 in the alphabet $\{1,2\}$, the one with smallest Markov value and with intruder sets we could found is
\begin{equation*}
    m\left(\overline{11221122^*1221111}\right)=3.09090912\dots,
\end{equation*}
which has period 15.

As a consequence of our explorations, we have confirmed that thousands of non semisymmetric odd orbits satisfy local uniqueness and self-replication, which implies the existence of thousands of regions of $M\setminus L$ with small dimension. More specifically, the result of our computer exploration is the following:

\begin{theorem}\label{thm:localself}
Let $\overline{w}\in(\N^*)^\Z$ be some periodic orbit where the minimal period $w$ is an odd non semisymmetric word of odd length less or equal than 39 and with Markov value $m(\overline{w})\leq\mu_2=3.037311\dots$. Then $w$ satisfies local uniqueness and self--replication, except for the following orbits and transposes of them:
\begin{gather*}
    211111221122111122211221111222112,
1111111222222211112221122211222211111, \\
122222112221122111122222112211112211221, 
122211222221122111122222112211221111221, \\
221111222221122111122222112211221122211,
211122112211221111122221122111112222112, \\
111111222112211111122211221111111221122,
112222112211221111122221122111112211122, \\
112211222211221111122221122111221111122, 
\end{gather*}
\end{theorem}

It is impossible for the exceptional orbits in \Cref{thm:localself} to satisfy local uniqueness (\Cref{def:local_uniqueness}). This failure arises for two distinct reasons: first, if a non semisymmetric minimal period $w$ of odd length attains its Markov value in more than one position, then local uniqueness cannot hold; second, if there exists two distinct pairs $\overline{u},\overline{u^T},\overline{v},\overline{v^T}$ of non semisymmetric orbits with the same Markov value, it can not hold either. In \Cref{sec:special_orbits} we classified all exceptional orbits from \Cref{thm:localself} in these two cases. These provide the first known examples of orbits exhibiting either of these phenomena. Remarkably, such orbits still produce elements of $M\setminus L$ and the corresponding regions remain free of intruder sets.

On the other hand, for the rest of the orbits, determining whether intruder sets exist requires substantially more computation, since it involves calculating the self-replication height. We have examined all previous orbits up to period 27 and found that each produces regions of $M\setminus L$ that are free of intruders, in the classical sense illustrated in \Cref{fig:MminusLclassic}. More precisely:

\begin{theorem}
Let $\overline{w}\in(\N^*)^\Z$ be some periodic orbit where the minimal period $w$ is an odd non semisymmetric word of odd length less or equal than 27 and with Markov value $m(\overline{w})\leq\mu_2=3.037311\dots$. Then $w$ satisfies local uniqueness and self--replication. Moreover, if we denote by $j_1=j_1(w)$ the self-replication height of $w$ (see \Cref{def:self-replication-height}), we have that
\begin{equation*}
    (M\setminus L)=C\cup D,
\end{equation*}
where $C$ is a topologically mixing Cantor set and $D$ is an infinite discrete set of points such that $D^\prime=C$. Moreover, we have the characterization
\begin{equation*}
    C = \left\{m(\underline{a})=\lambda_0(\underline{a})\in(j_0,j_1):\underline{a}\in(\N^*)^\Z, \min\{\ell(\underline{a}),\ell(\underline{a}^T)\}<j_0=\max\{\ell(\underline{a}),\ell(\underline{a}^T)\}\right\},
\end{equation*}
\begin{equation*}
    D = \left\{m(\underline{a})=\lambda_0(\underline{a})\in(j_0,j_1):\underline{a}\in(\N^*)^\Z, \ell(\underline{a})=\ell(\underline{a}^T)=j_0\right\}.
\end{equation*}
\end{theorem}

\begin{table}[h!]
\centering
    \begin{tabular}{C{3cm}C{5cm}}
        \hline
        Minimal period & Number of non semisymmetric orbits before $\mu_2$ \\ 
        \hline
        7 and 9 & 0 \\
        11 & 1 \\
        13 & 4 \\
        15 & 13 \\
        17 & 32 \\
        19 & 82 \\
        21 & 201 \\
        23 & 479 \\
        25 & 1131 \\
        27 & 2630 \\
        29 & 6084 \\
        31 & 14020 \\
        33 & 32283 \\
        35 & 74307 \\
        37 & 171171 \\
        39 & 394804
    \end{tabular}
    \caption{Number of non semisymmetric orbits with odd length that have Markov value less or equal than $\mu_2$ up to period 39. They all satisfy local uniqueness and self-replication: all of them give regions of $M\setminus L$. We have confirmed that up to period 27 they do not contain intruder sets.}
    \label{tab:regions_small_dimension}
\end{table}

\begin{question}\label{ques:small_intruders}
    Suppose $\overline{w}$ is a non semisymmetric orbit with odd length minimal period $w$ and Markov value $m(\overline{w})\leq t_{1/3}$ that satisfies local uniqueness and self--replication but such that the associated $M\setminus L$ region contains intruder sets? 
\end{question}

In \Cref{ques:small_intruders} we exclude the case where there exist two distinct pairs $\overline{u},\overline{u^T}$ and $\overline{v},\overline{v^T}$ of non semisymmetric orbits with the same Markov value $m(\overline{u})=m(\overline{v})$. In such a situation, local uniqueness (\Cref{def:local_uniqueness}) clearly fails.

The \Cref{ques:main} is whether the set $P$ intersects the Lagrange spectrum $L$. In general, we can ``approximate" portions of the Lagrange spectrum by sums of regular Cantor sets $K_1+K_2$ while the set $P$ is another regular Cantor set $K_3$, so the question is whether one has $(K_1+K_2)\cap K_3\neq\varnothing$. A relevant result from the PhD thesis \cite{JGThesis} is the following: for typical regular Cantor sets $K_1,K_2,K_3\subset\R$ with $\dim_H(K_1)+\dim_H(K_2)+\dim_H(K_3)>1$, if $\pi:\R^3\to\R$ is the affine projection $\pi(x,y,z)=x+y-z$, then $\pi(K_1\times K_2\times K_3)$ persistently contains interior.  In particular, for such typical triples of Cantor sets, the set of $t\in\R$ such that $(K_1+K_2)\cap (K_3+t)\neq\varnothing$ has non-empty interior. In the setting of \Cref{ques:main}, one has that each $\dim_H(K_i)$, $i=1,2,3$ is quite similar to $D(t)$, so one has that $\dim_H(K_1)+\dim_H(K_2)+\dim_H(K_3)>3\cdot D(t_{1/3})=1$, suggesting that an intersection between $P$ and $L$ may indeed occur. On the other hand, if one finds examples of regions in $M\setminus L$ with intruder sets corresponding to values of $t$ for which $D(t)<1/3$, then it may be possible to ``separate" the set $P$ from the Lagrange spectrum $L$.


\section{Definitions and notations}

A word $w$ of length $m$ is a finite string of positive integers $w=w_1\dots w_m\in\{1,\dots,A\}^m$. We denote the transpose word as $w^T:=w_m\dots w_1$. We say $w$ is palindrome if it coincides with its transpose $w=w^T$. We say that a periodic orbit $\overline{w}\in\{1,\dots,A\}^\Z$ has minimal period $w=w_1\dots w_m$ if $\overline{w}\neq\overline{w_1\dots w_n}$ for all $1\leq n<m$.

\begin{definition}\label{def:nonsemisymmetric}
    We say that $w\in\{1,\dots,A\}^{m}$ is \emph{semisymmetric} if it is a palindrome or the concatenation of two palindromes. Otherwise, it is \emph{non semisymmetric}.
\end{definition}

\begin{definition}
    We say that a periodic orbit $\overline{w}\in\{1,\dots,A\}^\Z$ is \emph{semisymmetric} if the orbits of $\overline{w}$ and $\overline{w^T}$ by the shift map are the same.
\end{definition}

It is a simple exercise to prove that a periodic orbit $\overline{w}\in\{1,\dots,A\}^\Z$ with period $w=w_1\dots w_m$ is semisymmetric if and only if $w$ is a semisymmetric word. It is also easy to see that most words are non semisymmetric. Through this paper we will only use $A=2$ or $A=3$.

Let $w=w_1\dots w_{2n+1}$ be a non semisymmetric word of odd length and denote $w^*=w_1\dots w_{n+1}^*\dots w_{2n+1}$. Suppose that this word is such that the Markov value is attained exactly at the middle position \[m(\overline{w})=\lambda_0(\overline{w_1\dots w_n w_{n+1}^{*}w_{n+2}\dots w_{2n+1}}).\]

\begin{definition}[Local uniqueness]\label{def:local_uniqueness}
    We say that $w\in\{1,\dots,A\}^{2n+1}$ satisfies local-uniqueness if there is an $\varepsilon_1(w)=\varepsilon_1>0$ such that for any bi-infinite sequence $\underline{\theta}\in\{1,\dots,A\}^{\Z}$ with $m(\underline{\theta})=\lambda_0(\underline{\theta})$, it holds that if $|m(\underline{\theta})-m(\overline{w})|<\varepsilon_1$, then $\underline{\theta}$ must be of the following form
    \begin{align*}
        \underline{\theta}&=\dots w_{n+1}\dots w_{2n+1} w_1 \dots w_{n}w_{n+1}^*w_{n+2}\dots w_{2n+1} w_1 \dots w_{n+1}\dots \\
        &= \dots  w_{n+1}\dots w_{2n+1} w^* w_1 \dots w_{n+1}\dots 
    \end{align*} 
    
\end{definition}

\begin{definition}[Self-replication]\label{def:self-replication}
    We say that $w\in\{1,\dots,A\}^{2n+1}$ satisfies self-replication to the left if there is an $\varepsilon_2(w)=\varepsilon_2>0$ such that for any bi-infinite sequence $\underline{\theta}\in\{1,\dots,A\}^{\Z}$ with $m(\underline{\theta})=\lambda_0(\underline{\theta})$, it holds that if $m(\underline{\theta})<m(\overline{w})+\varepsilon_2$ and $\underline{\theta}=\dots  w_{n+1}\dots w_{2n+1} w^* w_1 \dots w_{n+1}\dots$, then $\underline{\theta}$ must be of the following form
    \begin{align*}
        \underline{\theta}&=\dots w_{n+1}\dots w_{2n+1} w_1 \dots w_{n}w_{n+1}^*w_{n+2}\dots w_{2n+1} w_1 \dots w_{m}\dots \\
        &= \dots  w_{n+1}\dots w_{2n+1} ww^* w_1 \dots w_{m}\dots 
    \end{align*}
    for some $m<|w|+\lfloor |w|/2\rfloor=3n+1$ (with the convention $w_{k}:=w_{k-2n-1}$ when $k>2n+1$). In particular, for all such $\underline{\theta}$ it holds
    \begin{equation*}
        \underline{\theta}=\overline{w}w^*w_1\dots w_m\dots
    \end{equation*}
    
    We say $w$ satisfies self-replication to the right if the transpose word $w^T:=w_{2n+1}\dots w_1$ satisfies self-replication to the left.
\end{definition}

\begin{remark}
    It is clear that from the definition of self-replication (to some side) that one needs to repeat the period at least once to that side. The same is not true for the other side: it is possible to have self-replication even with $m<2n+1$. An explicit example of this phenomena occurs for the word $w=12111333112$ of length $n=11$, where $w$ self-replicates to the left with $m=2n$ and there are sequences $\underline{\theta}\in\{1,2,3\}^\Z$ with $m(\underline{\theta})<m(\overline{w})+\varepsilon_2(w)$ but $\underline{\theta}=\overline{w}w^*ww_1\dots w_{2n}3\dots$.  
\end{remark}

All the classical regions presented in \Cref{tab:knownregions} satisfy \Cref{def:local_uniqueness} and \Cref{def:self-replication}. For all these regions, it is possible to chose $\varepsilon_2(w)<\varepsilon_1(w)$. In fact, observe that this inequality guarantees that $L\cap(m(\overline{w})-\varepsilon_1(w),m(\overline{w})+\varepsilon_2(w))=\varnothing$. Indeed, since by \cite[Theorem 2, Chapter 3]{Cusick-Flahive} the Lagrange spectrum $L$ coincides with the closure of the Markov values of periodic points:
\begin{equation}\label{eq:lagrange_closure}
    L=\cl\left(\left\{m(\overline{a}):\overline{a}\in\Sigma \text{ completely periodic with period $a\in(\N^*)^n$}\right\}\right),
\end{equation}
if $\overline{p}$ is a periodic point with $|m(\overline{p})-m(\overline{w})|<\varepsilon_2(w)$, then since $\varepsilon_2(w)<\varepsilon_1(w)$ one has that $\overline{p}=\overline{w}w^*w_1\dots w_m\dots$ which can only occur if $\overline{p}=\overline{w}$.

For all these classical regions, the optimal value of $\varepsilon_2(w)$ always coincide with $j_1-j_0$, where $j_1=j_1(w)$ corresponds to the right border of the gap of $L$ to the right of $j_0$ (see \Cref{fig:MminusLclassic}). This can be seen because of the fact that the symbolic expression of $j_1$ in all these regions connects to the left and right with semi-symmetric orbits (thus different from $\overline{w}$), hence, the sequence $\underline{\theta}$ with $m(\underline{\theta})=j_1$ does not self-replicate $w$. For example, for $w=2112122$ one has $j_1=m\left(\overline{211212221211}ww^*ww\overline{122212112112}\right)$.

\begin{definition}[Self-replication height]\label{def:self-replication-height}
    Let $w\in\{1,\dots,A\}^{2n+1}$ be a finite word that satisfies self-replication to some side. Take the supremum $\varepsilon_2^{\max}(w)$ of all $\varepsilon_2(w)$ such that $w$ satisfies \Cref{def:self-replication}. The quantity $j_1=j_1(w):=m(\overline{w})+\varepsilon_2^{\max}(w)$ is called the \emph{self-replication height} of $w$.
\end{definition}

\begin{remark}
Notice that the self-replication height $j_1(w)$ depends on the orbit. If there are different pairs of non semisymmetric orbits $\overline{u},\overline{u^T}$ and $\overline{v},\overline{v^T}$ with the same Markov values $m(\overline{u})=m(\overline{v})$, their self-replication heights $j_1(u)$ and $j_1(v)$ are likely different. 
\end{remark}

For the new regions presented in this paper, one is not able to chose the parameter $\varepsilon_2(w)$ for self-replication satisfying $\varepsilon_2(w)<\varepsilon_1(w)$. For example, for $w=11122322221$, according to \Cref{fig:MminusLintruder1}, the number $i_0$ is the Markov value of a bi-infinite sequence that do not connect with $\overline{w}$ (see \Cref{lem1:branch_2}) and is such that $M\cap(j_0,i_0)=\varnothing$, so any such $\varepsilon_1(w)$ satisfies $\varepsilon_1(w)\leq j_0-i_0$. However, it is possible to prove that self-replication (to the right) holds for this particular word (\Cref{lem1:self-replication}), by choosing $\varepsilon_2(w)=j_1-j_0$, where is $j_1$ given in \Cref{lem:self_replication_words}. Again, the symbolic expression of the bi-infinite sequence with Markov value $j_1$ shows that this $\varepsilon_2(w)$ is optimal.

One of the goals of this paper is to show that even when the self-replication height is greater than the value needed for local uniqueness, i.e., when $m(\overline{w})+\varepsilon_1(w)<j_1(w)$, it is still possible to construct elements of $M\setminus L$. However, in contrast to the classical regions of \Cref{tab:knownregions}, it is very difficult to characterize explicitly all the elements of $(M\setminus L)\cap(j_0,j_1)$ that connect to the orbit $\overline{w}$.

The natural candidates for being in $M\setminus L$ are therefore given by
\begin{equation*}
    P:=\{m(B)=\lambda_0(B)\in(j_0,j_1): B=\overline{w}w^*w_1\dots w_m\ldots\in\{1,\dots,A\}^\Z\}.
\end{equation*}

For the classical regions in \Cref{tab:knownregions}, it is indeed the case that $(M\setminus L)\cap(j_0,j_1)=P$. However, for regions that exhibit intruder sets, one only expects that $P\subset(M\setminus L)\cap(j_0,j_1)$ because the existence of points in $L^\prime\cap(j_0,j_1)$ can produce regions of $M\setminus L$ associated with different non semi-symmetric orbits.

For most cases is very difficult to decide if one actually has $P\subset(M\setminus L)\cap(j_0,j_1)$. For the region associated with $w=11122322221$ that we study in detail in \Cref{sec:onegap}, one has that indeed $P\subset(M\setminus L)\cap(j_0,j_1)$, in fact $(M\setminus L)\cap[m_0,m_1]=P$. However, for the region associated with $w=1112122$, is very difficult to give a precise characterization. We are able to show for example that 
\begin{equation*}
    P\cap\left((b_1,a_2)\cup(b_2,a_3)\right)=(M\setminus L)\cap\left((b_1,a_2)\cup(b_2,a_3)\right)
\end{equation*}

A classical criteria to decide whether a Markov value belongs to $L^\prime$ was given by Flahive \cite[Theorem 2]{Flahive77}. 

\begin{lemma}\label{lem:flahive}
    Let $\underline{a}\in(\N^*)^\Z$ be a bi-infinite sequence such that $m(\underline{a})=\lambda_0(\underline{a})$. Suppose that $\underline{a}$ is eventually periodic at both sides with semi-symmetric periods, that is, $\underline{a}=\overline{p}~a_{-m}\dots a_{n}~\overline{q}$ with $p,q$ semi-symmetric words. If $m(\underline{a})\neq m(\overline{p}),m(\overline{q})$, then $m(\underline{a})\in L^\prime$.
\end{lemma}

\begin{lemma}
Let $w\in\{1,\dots,A\}^{2n+1}$ be a word that satisfies self-replication with self-replication height $j_1$. If $t$ is such that $m(\overline{w})\leq t<j_1$, then $\Sigma_t$ is not transitive.
\end{lemma}

\begin{proof}
Since all periodic orbits $\overline{w}$ with Markov value $m(\overline{w})\leq 3$ have even length and are semisymmetric \cite{Bombieri}, assume that $m(\overline{w})>3$. Suppose there is a bi-infinite sequence $\underline{\theta}=(\theta_n)_{n\in\Z}$ such that the closure of $\{\sigma^n(\underline{\theta}):n\in\Z\}$ is equal to $\Sigma_t$. Since $\overline{w}\in\Sigma_t$, there is $k\in\Z$ and $\ell\in\Z$ such that the block $\theta_{k-2n}\dots\theta_{k+2n}$ coincides with $w_{n+2}\dots w_{2n+1}w_1\dots w_{2n+1}w_1\dots w_n$ and the block $\theta_{\ell-2n}\dots\theta_{\ell+2n}$ coincides with the transpose, that is, coincides with $w_{n}\dots w_{1}w_{2n+1}\dots w_{1}w_{2n+1}\dots w_{n+2}$. From \Cref{def:self-replication} we obtain that $\sigma^{k}(\underline{\theta})$ and $\sigma^\ell(\underline{\theta})$ must have the form
\begin{equation*}
    \sigma^k(\underline{\theta})=\overline{w}w^*w_1\dots w_m\dots, \quad\text{and}\quad \sigma^\ell(\underline{\theta})=\dots w_m\dots w_1w^*\overline{w^T},
\end{equation*}
for some positive integer $m$, respectively. Since $\overline{w}$ and $\overline{w^T}$ are distinct, we must have that $k<\ell$ and so $\underline{\theta}$ has the form $\overline{w}\tau\overline{w^T}$ for some finite word $\tau$. Since the closure of such an orbit is countable and $\Sigma_t$ is uncountable for $t\geq 3$, we have a contradiction, so we conclude that $\Sigma_t$ is not transitive.
\end{proof}


\section{$M\setminus L$ associated with a simple intruder set}\label{sec:onegap}

Through this section $w:=11122322221$ and $w^*:=111223^*22221$. We have that $j_0:=m(\overline{111223^*22221})=3.8340731702358\dots$.

We will show that there is a region of $M\setminus L$ that connects with $w$ and which has the structure of \Cref{fig:MminusLintruder1}.

\subsection{Local uniqueness and self-replication}

\begin{lemma}[Local uniqueness forbidden words]\label{lem1:local_uniqueness_words}
    The following words 
    \begin{gather*}
        13, 32323, 32322, 1223221, 21223222, 1112232223, 21112232222, \\
        11w, w2, 321w11, 221w112,  1221w11, 221w1111,  22221w11123, 32221w11123, \\
        32221w11122, 332221w11121121, 232221w1112112 \\
        2211111223222222,
        1211111223222222,
        33211111223222222, \\
        211111223222223, 2111112232222221, 2111112232222222, 21111122322222233 \\
        123211111223222222321, 2123211111223222222322, 
    \end{gather*}
    produce Markov values at least $3.8340731707>j_0+10^{-10}$.
\end{lemma}

\begin{lemma}[Local uniqueness]\label{lem1:local_uniqueness}
    Suppose $B\in\{1,2,3\}^\Z$ satisfies $m(B)<3.8340731707$ and $3.834073167<\lambda_0(B)<3.8340731707$. Then $B=22221w^*11122$ or $B=12221w^*11123$.
\end{lemma}

\begin{proof}
    Since $\lambda_0(33^*3)<3.62$, $\lambda_0(23^*3)<3.75$ and 13 is forbidden, we can assume $B=23^*2$. Since $\lambda_0(123^*2)<3.833$ we must continue either as $223^*22$ or $323^*21$ (since 32322, 32323 are forbidden). Thus we can assume that $B=223^*22$. Since $\lambda_0(2223^*222)<3.831$, $\lambda_0(3223^*223)<3.822$, $\lambda_0(3223^*2211)<3.832$, $\lambda_0(3223^*2212)<3.834$ and 13, 1223221 are forbidden, we can assume that $B=1223^*222$. Since 21223222 and 31 are forbidden we must extend to $B=11223^*222$. If it continues as $211223^*222$, then using 31 is forbidden the value would be less than $\lambda_0(211223^*22232)<3.834$. Thus $B=111223^*222$. Since $\lambda_0(111223^*2221)<3.834$ and 1112232223, 21112232223, 31 are forbidden, we must extend to $B=1111223^*2222$. We have then two cases:
    \begin{itemize}
        \item If it continues as $B=1111223^*22221=1w^*$, then using the forbidden words $31,11w,w2$ it will extend to $B=21w^*1$. If it continues as $21w^*12$, then we will have that $\lambda_0(B)\leq \lambda_0(\overline{32}1w^*1\overline{23})<3.83407$. Thus $B=21w^*11$. If it continues as $121w^*11$, then since $\lambda_0(B)\leq\lambda_0(\overline{21}121^*11\overline{12})<3.83406$ and $321w11$ is forbidden, we must continue as $B=221w^*11$. Since $13, 221w112, 221w1111, 1221w11$ are forbidden and $\lambda_0(3221w^*1112)<3.8340717$ we must continue as $B=2221w^*1112$. Since $32221w11123, 22221w11123,3222w11122$ are forbidden and $\lambda_0(12221w^*11121)<3.8340724$, $\lambda_0(22221w^*11121)<3.8340731$, $\lambda_0(12221w^*11122)<3.8340731$ we must continue either as $B=22221w^*11122$ or $B=12221w^*11123$ or $B=32221w^*11121$. In the particular case when $B=32221w^*11121$, since $\lambda_0(32221w^*111212)<3.834073$, $\lambda_0(332221w^*1112111)<3.8340731$, $\lambda_0(332221w^*11121123)<3.83407316$, $\lambda_0(332221w^*11121122)<3.834073167$, $\lambda_0(2332221w^*11121122)<3.834073167$, $\lambda_0(232221w^*1112111)<3.83407314$ and $332221w11121121, 232221w1112112,13$ are forbidden, we see that this case has no extensions. 
        
        \item If it continues as $B=1111223^*22222$, then since \[\lambda_0(21111223^*22222)<3.834,\] \[\lambda_0(111111223^*22222)<3.83407314\] and 31 is forbidden, we must continue as $B=211111223^*22222$. Since 211111223222223 and \[\lambda_0(21211111223^*22222112)<3.834072,\] we must continue as $B=211111223^*222222$. Now since 2111112232222221, 2111112232222222, 21111122322222233, 2211111223222222, 1211111223222222, 33211111223222222, 13 are forbidden it must extend to $B=23211111223^*22222232$. Since 
        \begin{align*}
            \lambda_0(323211111223^*22222232)&<3.8340731, \\
            \lambda_0(223211111223^*22222232)&<3.83407317,
        \end{align*}
        we must extend to $B=123211111223^*22222232$. Since 123211111223222222321 and 22323 are forbidden we must continue as $B=123211111223^*222222322$. Since \[\lambda_0(1123211111223^*222222322)<3.834073166,\] and 2123211111223222222322, 31 are forbidden, this word has no extensions.
    \end{itemize}
\end{proof}

\begin{lemma}[Self-replication forbidden words]\label{lem:self_replication_words}
    The following words
    \begin{gather*}
        122221w11122, 222221w11122, 22221w111221, 22221w111222, 2322221w1112233, \\ 
        3322221w1112232, 3322221w1112233, 12322221w1112232, 2322221w11122321, \\
        222322221w111223222, 222322221w111223223,  322322221w111223222, \\
        322322221w111223223, 
        1122322221w1112232221, 1122322221w111223223, \\
        21122322221w1112232222, 1ww11122322223, 21ww11122322222, 121www111, \\
        121www12, 221www12, 3221www1112, 2221www11121
    \end{gather*}
    give Markov values at least 
    \begin{align*}
        j_1 := m(\overline{21}1232321212112221www^*w11121\overline{12})&=3.83407317023586722\ldots \\
        &\approx j_0+2.3728\cdot10^{-14},
    \end{align*}
    which corresponds to the minimum Markov value of a bi-infinite sequence containing $2221www11121$.

\end{lemma}

\begin{lemma}[Self-replication]\label{lem1:self-replication}
    Suppose $B\in\{1,2,3\}^\Z$ satisfies $B=22221w^*11122$ and $m(B)<j_1$. Then $B=2221ww^*w11122$. 
\end{lemma}

\begin{corollary}\label{cor1:first_characterization}
    Suppose $B\in\{1,2,3\}^\Z$ satisfies $B=22221w^*11122$ and $m(B)<j_1$. Then $B=2221ww^*\overline{w}$.
\end{corollary}

\begin{lemma}\label{lem:good_characterization}
    Suppose $B\in\{1,2,3\}^\Z$ satisfies $B=22221w^*11122$ and $ j_0<m(B)=\lambda_0(B)<j_1$. Then  $B=21212112221ww^*\overline{w}$ and $m_0\leq \lambda_0(B)\leq m_1$, where
    \begin{align*}
    m_0 := m(\overline{21}12221ww^*\overline{w})    &=3.83407317023586720605935\ldots \\
    &\approx j_0+2.3708\cdot10^{-14},
\end{align*}
\begin{align*}
    m_1:=m(\overline{21}1232321212112221ww^*\overline{w})&= 3.834073170235867216\ldots \\
    &\approx j_0+2.3718\cdot10^{-14}.
\end{align*}

\end{lemma}

\begin{proof}

From \Cref{cor1:first_characterization} we know that $B=2221ww^*\overline{w}$. If $B$ does not contain $12221www11122$, then using inductively \Cref{lem1:self-replication} we will get that $B=\overline{w}$, a contradiction. So we must have that $B=12221ww^*\overline{w}$ and thus $m(B)\geq m_0$, where we used that 13, 31 are forbidden. Since
    \begin{align*}
        \lambda_0(212221ww^*\overline{w})>j_1+10^{-14}, \\
        \lambda_0(1112221ww^*\overline{w})>j_1+10^{-15}, \\
        \lambda_0(32112221ww^*\overline{w})>j_1+10^{-15}, \\
        \lambda_0(22112221ww^*\overline{w})>j_1+10^{-15}, \\
        \lambda_0(112112221ww^*\overline{w})>j_1+10^{-16}, \\
        \lambda_0(3212112221ww^*\overline{w})>j_1+10^{-16}, \\
        \lambda_0(2212112221ww^*\overline{w})>j_1+10^{-17}, \\
        \lambda_0(11212112221ww^*\overline{w})>j_1+10^{-18},
    \end{align*}
    and 13 is forbidden, we see that we must continue as $B=21212112221ww^*\overline{w}$.
    
    Finally, use that if 13, 32322, 32323 are forbidden, then a continued fraction in $\{1,2,3\}$ that begins with $[0;2]$ is maximized with $[0;2,3,2,3,2,1,\overline{1,2}]$.
\end{proof}

Let us denote our candidate set for points of $M\setminus L$ by
\begin{equation*}
	P:=\{m(B)=\lambda_0(B)\in(j_0,j_1): B=\dots21212112221ww^*\overline{w}\in\{1,2,3\}^\Z\}.
\end{equation*}

Notice that $m_0=\min P$ and $m_1=\max P$. We will show that $P\subset M\setminus L$. 

\subsection{Separating the intruder combinatorics}

Let us call $\tau:=12221w11123$ and $\tau^*:=12221w^*11123$. This continuation correspond to the \emph{intruder combinatorics}. We would like to show that $\lambda_0(\tau^*)$ is different from $\lambda_0(22221w^*11122)$ for any allowed continuations with values below some threshold. This upper bound is naturally given by the self-replication height. 

\begin{lemma}\label{lem1:second_branch_forbidden_words}
    The following words
    \begin{gather*}
        221\tau2,
        121\tau2,
        21\tau22,
        21\tau23,
        2121\tau, \\
        21\tau211,
        11121\tau,
        3321\tau2,
        321121\tau, 
        221121\tau,\\
        121\tau33,
        121\tau321,
        121\tau322,
        21\tau2123,
        21\tau2122, \\
        212321\tau21211,
        12321\tau212112,
        21112321\tau2121111, \\
        1112321\tau21211112,
        2111112321\tau212111111,
        111112321\tau2121111112, \\
        11111112321\tau212111111112,
        211111112321\tau21211111111, \\
        1111111112321\tau21211111111,
        11111112321\tau2121111111111
    \end{gather*}
    give Markov values at least 
    \begin{equation*}
        m(\overline{21}1111111112321\tau^*21211111111\overline{12})=3.834073170235887\ldots\approx j_0+4.39\cdot10^{-14}
    \end{equation*}
    which corresponds to the minimum Markov value of a bi-infinite sequence that contains $1111111112321\tau21211111111$.
\end{lemma}

\begin{lemma}\label{lem1:branch_1}
    Suppose $B\in\{1,2,3\}^\Z$ satisfies $B=\tau^*$. Then $|m(B)-j_0|>4\cdot 10^{-14}$ or $B=2111111112321\tau^*2121111111112$.
\end{lemma}
\begin{proof}
    Since $\lambda_0(2\tau^*)<3.834073$ and 31 is forbidden we must continue as $B=1\tau^*$. Since $\lambda_0(11\tau^*)$ is at most $\lambda_0(\overline{21}11\tau^*\overline{23})<3.83407314$, we can assume it continues as $B=21\tau^*$. We have two cases:
    \begin{itemize}
        \item If we continue as $B=21\tau^*3$, then since $\lambda_0(321\tau^*3)<3.83407316$, $\lambda_0(221\tau^*3)<3.83407317$ it must continue as $B=121\tau^*3$. Now since $2121\tau$, $11121\tau$, $221121\tau$, $321121\tau$, $121\tau33$, $121121\tau321$, $121121\tau322$, $31$ are all forbidden this is forced to $B=121121\tau^*323$, which satisfies
        \begin{equation*}
            \lambda_0(B)\leq m(\overline{21}121121\tau^*32332321\overline{12})=3.83407317023006\dots
        \end{equation*}
        \item If we continue as $B=21\tau^*2$, since $221\tau2$, $121\tau2$,
        $3321\tau2$,
        $21\tau22$, $21\tau23$, $21\tau211$, $21\tau2122$, $21\tau2123$, $13$ are all forbidden we are forced to $B=2321\tau^*2121$. Now since 
        \begin{align*}
            \lambda_0(32321\tau^*2121)&<3.83407317002, \\
            \lambda_0(22321\tau^*2121)&<3.83407317021, \\
            \lambda_0(12321\tau^*21212)&<3.8340731702,
        \end{align*}
        and $13$ is forbidden we must extend to  $B=12321\tau^*21211$. Now since $212321\tau21211$, $12321\tau212112$, $31$ are forbidden and 
        \begin{align*}
            \lambda_0(2112321\tau^*212111)&<3.8340731702317, \\
            \lambda_0(1112321\tau^*2121112)&<3.83407317023,
        \end{align*}
        we must continue as $B=1112321\tau^*2121111$. Moreover since $21112321\tau2121111$, $1112321\tau21211112$, $31$ are forbidden and 
        \begin{align*}
            \lambda_0(211112321\tau^*21211111)&<3.8340731702355, \\
            \lambda_0(11112321\tau^*212111112)&<3.8340731702352,
        \end{align*}
        we must continue as $B=111112321\tau^*212111111$. Now since \[2111112321\tau212111111, 111112321\tau2121111112, 13\] are forbidden we are forced to $B=1111112321\tau^*2121111111$. Since
        \begin{align*}
            \lambda_0(21111112321\tau^*212111111112)&<3.83407317023579, \\
            \lambda_0(211111112321\tau^*21211111112)&<3.83407317023575
        \end{align*}
        we must continue as $B=11111112321\tau^*21211111111$. Finally, since 
        \begin{gather*}
            211111112321\tau21211111111, 11111112321\tau212111111112,  \\
            1111111112321\tau21211111111,
        11111112321\tau2121111111111, 31
        \end{gather*} 
        are forbidden we must continue as $B=2111111112321\tau^*2121111111112$.   
    \end{itemize}
\end{proof}

\begin{lemma}\label{lem1:branch_2}
    Suppose $B\in\{1,2,3\}^\Z$ satisfies $B=2111111112321\tau^*2121111111112$ and $m(B)=\lambda_0(B)$. Then either
    \begin{align*}
        m(B)&\leq m(\overline{12}2111111112321\tau^*212111111111232321\overline{12})\\
        &=3.8340731702358417\ldots\approx j_0-1.75\cdot 10^{-15}
    \end{align*}
    or $i_0\leq m(B)\leq i_1$ where
    \begin{align*}
        i_0 &:= m(\overline{21}1232332111111112321\tau^*2121111111112332321\overline{12}) \\
        &=3.834073170235847\ldots\approx j_0+4.08\cdot10^{-15},
    \end{align*}
    and
    \begin{align*}
        i_1 &:= m(\overline{21}123232111111112321\tau^*212111111111232321\overline{12}) \\
        &=3.834073170235850\ldots\approx j_0+7.28\cdot10^{-15}.
    \end{align*}
    
\end{lemma}

\begin{proof}
    We will use that if 13, 32322, 32323 are forbidden, then a continued fraction in $\{1,2,3\}$ that begins with $[0;2]$ is maximized with $[0;2,3,2,3,2,1,\overline{1,2}]$ and similarly a continued fraction in $\{1,2,3\}$ that begins with $[0;2,3]$ is minimized with $[0;2,3,3,2,3,2,1,\overline{1,2}]$.
    
    Since \[\lambda_0(12111111112321\tau^*2121111111112)<3.834073170235825,\] \[\lambda_0(2111111112321\tau^*21211111111121)<3.834073170235822,\]
    we see that have that either $B=22111111112321\tau^*21211111111123$ and so
    \begin{equation*}
        \lambda_0(B)\leq m(\overline{12}2111111112321\tau^*212111111111232321\overline{12})=3.8340731702358417\ldots
    \end{equation*}
    or that either $B=32111111112321\tau^*21211111111122$
    \begin{equation*}
        \lambda_0(B)\leq m(\overline{21}123232111111112321\tau^*2121111111112\overline{21})=3.8340731702358414\ldots
    \end{equation*}
    or that either  $B=32111111112321\tau^*21211111111123$, for which $m(B)$ has values bounded between $i_0$ and $i_1$.
\end{proof}

In particular we obtain
\begin{lemma}\label{lem1:second_characterization}
    Suppose $B\in\{1,2,3\}^\Z$ satisfies $i_1< \lambda_0(B)\leq m(B)<j_1$. Then up to transposition $B=2221ww^*\overline{w}$ and if $m_0\leq \lambda_0(B)$ then $B=21212112221ww^*\overline{w}$ and
    \begin{align*}
        m_0\leq \lambda_0(B)\leq m_1.
    \end{align*}
\end{lemma}
\begin{proof}
    By \Cref{lem1:local_uniqueness}, we have that either $B=22221w^*11122$ or $B=\tau^*$. In the latter case, by \Cref{lem1:branch_1} and \Cref{lem1:branch_2} we must have that either $m(B)\leq i_1<m_0$ or that $m(B)<j_0-10^{-15}$, which are impossible. Thus $B=22221w^*11122$, so \Cref{lem:good_characterization} finishes the proof.
\end{proof}

\begin{corollary}\label{cor1:gaps}
    The interval $(i_1,j_1)$ is a maximal gap of $L$ and $(i_1,m_0)$ is a maximal gap of $M$.
\end{corollary}
\begin{proof}
    Suppose $B\in\{1,2,3\}^\Z$ and $m(B)=\lambda_0(B)\in (i_1,j_1)$. By \Cref{lem1:second_characterization}, we can assume that $B=2221ww^*\overline{w}$. If $B$ is periodic then $B=\overline{w^*}$, a contradiction. In case that $\lambda_0(B)<m_0$, then $B$ can not contain the word $12221www11122$, so applying inductively \Cref{lem1:self-replication} we will also get $B=\overline{w^*}$. 
\end{proof}

\begin{corollary}
    Suppose $B\in\{1,2,3\}^\Z$ satisfies $j_0-10^{-15}<\lambda_0(B)=m(B)<i_0$. Then $B=\overline{w^*}$ and $m(B)=j_0$.
\end{corollary}
\begin{proof}
    By \Cref{lem1:local_uniqueness}, we have that either $B=22221w^*11122$ or $B=\tau^*$. In the former case, by \Cref{cor1:first_characterization} we must have $B=2221ww^*\overline{w}$. Since $B$ can not contain $12221www11122$, applying inductively \Cref{lem1:self-replication} we see that $B=\overline{w^*}$. In the latter case, by \Cref{lem1:branch_1} and \Cref{lem1:branch_2} we must have that either $m(B)\geq i_0$ or that $m(B)<j_0-10^{-15}$, which are impossible.
\end{proof}

We can characterize points in $M\setminus L$ associated with $w$ in this region. Define $F$ to be the set of all words of \Cref{lem1:local_uniqueness_words}, \Cref{lem1:second_branch_forbidden_words} and the self-replicating word $22221w11122$, and their transposes. The proof of the following explicit characterization theorem is done in the same way as in the classical regions of \Cref{tab:knownregions}, but we include it for completeness.

\begin{theorem}\label{thm1:characterization}
    We have that
    \begin{equation*}
        (M\setminus L)\cap(i_1,j_1)=C\cup D\cup X
    \end{equation*}
    where 
    \begin{multline*}
        C=\{\lambda_0(\gamma21212112221ww^*\overline{w}):\gamma21212112221ww\in\{1,2,3\}^{\N} \\
        \text{does not contain any word from $F$}\},
    \end{multline*}
    is a Cantor set and
    \begin{multline*}
        D=\big\{\lambda_0(\overline{w^T}ww12221121212\theta21212112221ww^*\overline{w}):\theta~\text{finite word in 1,2,3, } \\ [0;\theta]\geq[0;\theta^T], \text{and } 12221121212\theta21212112221 \text{ not contains any word from $F$}\big\},
    \end{multline*}
    \begin{multline*}
        X=\big\{\lambda_0(\overline{w^T}ww122211212121212112221ww^*\overline{w}), \\
        \lambda_0(\overline{w^T}ww1222112121212112221ww^*\overline{w}), \lambda_0(\overline{w^T}ww12221121212112221ww^*\overline{w})\big\},
    \end{multline*}
    are sets of isolated points in $M$.
\end{theorem}
\begin{proof}
    By \Cref{lem1:second_characterization} and \Cref{cor1:gaps} we have that $B=\gamma21212112221ww^*\overline{w}$ for some $\gamma\in\{1,2,3\}^\N$. Note that $\gamma21212112221$ can not contain the self-replicating word $22221w11122$ because otherwise $12221$ would be contained in $\overline{w}$. In case it contains the transpose $22111w^T12222$, we can apply \Cref{lem1:second_characterization} to $B^T$ to see that either $\lambda_0(B)$ belongs to $D$ or $X$. Otherwise it must belong to $C$. To see the converse inclusion, for $\lambda_0(B)\in D\cup X$ is direct computation. For $\lambda_0(B)\in C$, note that if $i_1<\lambda_i(B)$ for some position inside $\gamma$, then by \Cref{lem1:second_characterization} we will get that $\gamma$ contains $22111w^T12222$ contrary to the definition of $C$.
\end{proof}



\section{$M\setminus L$ associated to an intruder set difficult to separate}\label{sec:length7}

Lets call $w:=1112122$ and $w^*:=1112^*122$. We have that $j_0:=m(\overline{1112^*122})=3.334156277015\dots$.

\subsection{Local uniqueness and self-replication}

\begin{lemma}[Local uniqueness forbidden words]\label{lem3:local_uniqueness_forbidden_words}
    All the words
    \begin{gather*}
        21212, 12121112, 12111212, 111112121, 111212111, 2w2,  2w112, 222w11,  \\ 2211112121122, 21121121211112212, 1212211112121121111, \\
        2122111121211211112, 221122111121211211212, 
    \end{gather*}
    give Markov values at least $3.3341579>j_0+10^{-6}$. 
\end{lemma}

\begin{lemma}[Local uniqueness]\label{lem3:local_uniqueness}
    Suppose $B\in\{1,2\}^\Z$ satisfies $3.33415618<\lambda_0(B)=m(B)<3.3341568$. Then $B=122w^*111$ or $B=22212211112^*121121111112$. 
\end{lemma}

\begin{proof}
    Since $\lambda_0(1^*)<2.4$, $\lambda_0(2^*2)<3.16$ we must have $B=12^*1$. Since 21212 is forbidden and $\lambda_0(112^*11)<3.27$, we must continue as $B=112^*12$. Since $\lambda_0(2112^*12)<3.32$, we must continue as $B=1112^*12$. We have two cases:
    \begin{itemize}
        \item If it continues as $1112^*122=w^*$, the using that $\lambda_0(1w^*)<3.329$ and $2w2$ is forbidden, we must continue as $2w^*1$. Since 12111212 is forbidden and $\lambda_0(22w^*12)<3.33397$, we must continue as $22w^*11$. Finally since $2w112$ and $222w11$ are forbidden we are forced to extend to $122w^*111$.
        \item If it continue as $1112^*121$, then using that 21212 and 21112121 are forbidden we must continue as $11112^*1211$. Since 111112121 and 111212111 are forbidden we must continue as $211112^*12112$. Since $\lambda_0(1211112^*12112)<3.3339$ and 2211112121122 is forbidden we must continue as $2211112^*121121$. Since $\lambda_0(2211112^*1211212)<3.3341$ and $\lambda_0(22211112^*121121)<3.3341$ we must continue as $12211112^*1211211$. Since $\lambda_0(112211112^*12112111)<3.33413$ and 21221111212112112 is forbidden, we must continue as either $112211112^*12112112$ or $212211112^*12112111$. We have two subcases:
        \begin{itemize}
            \item In the case $112211112^*12112112$, using that 111221111212112112 is forbidden and $\lambda_0(2112211112^*121121122)<3.334144$, it must continue as $2112211112^*121121121$. Since $\lambda_0(2112211112^*1211211211)<3.334156$ and 221122111121211211212 is forbidden, we must continue as $12112211112^*1211211212$. Finally since \[\lambda_0(12112211112^*12112112122)<3.334156,\] and 21212 is forbidden, we must continue as \[\lambda_0(12112211112^*121121121211)<3.33415618.\] 
            
            \item In the case $212211112^*12112111$, using that $\lambda_0(212211112^*121121112)<3.33415$ and 1212211112121121111 is forbidden we must continue as $2212211112^*121121111$. Since 2122111121211211112 is forbidden and $\lambda_0(12212211112^*1211211111)<3.3341555$, we must continue as $22212211112^*1211211111$. Since $\lambda_0(22212211112^*12112111112)<3.334153$, $\lambda_0(22212211112^*121121111111)<3.3341558$, we must continue as $22212211112^*121121111112$.

        \end{itemize}

    \end{itemize}
\end{proof}

\begin{lemma}[Self replication forbidden words]\label{lem3:self_replication_forbidden_words}
    All the words
    \begin{equation*}
        1122w11, 2w1111, 2w11122, 22122w11, 112122w111211, 2112122w11, 1www,
    \end{equation*}
    give Markov values at least
    \begin{align*}
        j_1:=m(\overline{212111112122}w^*ww121112\overline{221211111212})&=3.3341574813\dots\\
        &\approx j_0+1.204\cdot10^{-6}
    \end{align*}
    which corresponds to the minimum Markov value of a bi-infinite sequence containing $1www$.
\end{lemma}

\begin{lemma}[Self-replication]\label{lem3:self-replication}
    Suppose $B\in\{1,2\}^\Z$ satisfies $B=22w^*11$, $m(B)<j_1$. Then $B=22ww^*w1$. 
\end{lemma}

\begin{proof}
    Since $2w112$ and $222w11$ are forbidden, it must extend to $122w^*111$. Since $1122w11$ and $2w1111$ are forbidden, it must extend to $2122w^*1112$. Since $22122w11$ and $2w11122$ are forbidden, it must extend to $12122w^*11121$. Since $21212$ and $12122w111211$ are forbidden, it must extend to $112122w^*111212$. Since $2112122w11$ and $21112121$ are forbidden, it must extend to $ww^*w$. Since $1www$ and $2w2$ are forbidden, it must extend to $2ww^*w1$. Since 12111212 is forbidden, it must extend to $22ww^*w1$. 
\end{proof}

\begin{corollary}\label{cor3:first_characterization}
    Suppose $B\in\{1,2\}^\Z$ satisfies $B=22w^*11$, $m(B)<j_1$. Then $B=\overline{w}w^*w1$. Moreover if $j_0<m(B)=\lambda_0(B)<j_1$ then $B=\overline{w}w^*w12$.
\end{corollary}

\begin{proof}
    Applying inductively \Cref{lem3:self-replication} we must connect to the left with the periodic orbit $\overline{w}$, namely $B=\overline{w}w^*w1$. If the sequence continues as $B=\overline{w}w^*w11\dots$, then since $2w112$ is forbidden, applying \Cref{lem3:self-replication} inductively one is forced to continue as $B=\overline{w}w^*ww^k12\dots$ for some $k\geq 1$ (since $j_0\neq m(B)$), but this contradicts that $m(B)=\lambda_0(B)$ because in this situation $\lambda_0(\overline{w}w^*w12\dots)=\lambda_{7k}(B)>\lambda_0(B)$.
\end{proof}

Define 
\begin{equation*}
    m_0:= m(\overline{w}w^*w1221211\overline{122212111112})=3.334156382\dots\approx j_0+1.05\cdot10^{-7}.
\end{equation*}
which corresponds to the minimum Markov value of a word containing $22www12$ and also define
\begin{equation*}
    m_1:=m(\overline{w}w^*w1211\overline{122212111112})\approx 3.33415641441028\dots\approx j_0+1.37\cdot10^{-7}
\end{equation*}

We will see that the numbers $m_0$ and $m_1$ belong to $(M\setminus L)\cap(j_0,j_1)$. These numbers bound the values of the candidates that connect to $\overline{w}$.

\begin{lemma}\label{lem3:bounding_P_values}
    Suppose $B\in\{1,2\}^\Z$ satisfies $B=22w^*11$ and $j_0< \lambda_0(B)=m(B)<j_1$. Then $B=\overline{w}w^*w12$ and $m_0\leq m(B)\leq m_1$.
\end{lemma}

\begin{proof}
    Use the previous corollary and the forbidden words $2w2$, 111112121, 21212, 21211121.
\end{proof}

\subsection{Separating intruder sets from $M\setminus L$}

Let us denote
\begin{equation*}
    \tau:=22212211112121121111112,\quad\text{and}\quad \tau^*:=22212211112^*121121111112.
\end{equation*}
This continuation correspond to the \emph{intruder combinatorics}. We would like to show that the possible values of $\lambda_0(\tau^*)$ are different from $m_0$ and $m_1$. Note that this automatically gives that $m_0,m_1\in M\setminus L$. In fact, we will show that the distance of $\lambda_0(\tau^*)$ to $m_0$ and $m_1$ is positive, so we will be able to show the existence of Cantor sets of $M\setminus L$ containing $m_0$ and $m_1$ respectively.

More generally, we would like to show that
\begin{equation*}
    P:=\{m(B)=\lambda_0(B)\in(j_0,j_1): B=\overline{w}w^*w12\dots\in\{1,2\}^\Z\}
\end{equation*}
is contained in $(M\setminus L)\cap(j_0,j_1)$. For this, we must show that for any $\lambda_0(\tau*)\in(j_0,j_1)$ it holds that $\lambda_0(\tau^*)\not\in P$. However, this turns out to be very difficult to achieve. For our purposes, we will only separate the intruder combinatorics from the minimum and maximum of $P$. 

Notice that $m_0 = \min P$ and $m_1 = \max P$, so the set $P$ can be covered by the interval $[m_0,m_1]$. To separate $P$ from the intruder sets, we need to refine the covering. For that purpose, given a finite word $a_1\dots a_n\in\{1,2\}^\Z$ such that there is at least a sequence $S=\overline{w}w^*w12a_1\dots a_n\dots\in\{1,2\}^\Z$ with $m(S)=\lambda_0(S)\in(j_0,j_1)$, denote
\begin{align*}
    p_{a_1\dots a_n}^{(0)} &= \min\{m(B)=\lambda_0(B)\in(j_0,j_1):B=\overline{w}w^*w12a_1\dots a_n\dots\in\{1,2\}^\Z\}, \\
    p_{a_1\dots a_n}^{(1)} &= \max\{m(B)=\lambda_0(B)\in(j_0,j_1):B=\overline{w}w^*w12a_1\dots a_n\dots\in\{1,2\}^\Z\}
\end{align*}

Computing those Markov values reduces to minimize or maximize certain continued fractions that avoid all the local uniqueness forbidden words. 

\subsubsection{Separating $m_0$ from the intruder}

\begin{lemma}\label{lem3:intruder_branch_1}
    Suppose $B\in\{1,2\}^\Z$ satisfies, $B=\tau^*$, $m(B)\leq m_1$. Then $|m(B)-m_0|>7\cdot 10^{-11}$.
\end{lemma}

\begin{proof}
    Since $\lambda_0(1\tau^*)>m_0+10^{-7}$ and $\lambda_0(\tau^*1)>m_0+10^{-7}$ we must continue as $2\tau^*2$. Since $\lambda_0(22\tau^*2)>m_0+2.8\cdot10^{-8}$, $\lambda_0(12\tau^*22)>m_0+1.5\cdot10^{-8}$ we must continue as $12\tau^*21$.  Since $\lambda_0(212\tau^*21)<m_0-10^{-8}$ and $\lambda_0(12\tau^*212)<m_0-5\cdot 10^{-9}$ we must continue as $112\tau^*211$. Since $\lambda_0(2112\tau^*211)>m_0+4\cdot 10^{-9}$ and $\lambda_0(112\tau^*2112)>m_0+2\cdot 10^{-9}$ we must continue as $1112\tau^*2111$. Since $\lambda_0(1112\tau^*21112)<m_0-2\cdot 10^{-10}$ and $\lambda_0(21112\tau^*2111)<m_0-9\cdot 10^{-10}$ then we must continue as $11112\tau^*21111$. Since $\lambda_0(11112\tau^*211112)>m_0+10^{-9}$ and $\lambda_0(211112\tau^*21111)>m_0+10^{-9}$ then we must continue as $111112\tau^*211111$. Finally, since $\lambda_0(1111112\tau^*211111)>m_0+10^{-10}$, $\lambda_0(2111112\tau^*2111111)>m_0+7\cdot 10^{-11}$ we must continue as $\lambda_0(2111112\tau^*2111112)<m_0-5\cdot 10^{-11}$   
\end{proof}

\begin{lemma}\label{lem3:intruder_branch_2}
    Suppose $B\in\{1,2\}^\Z$ satisfies $B=\tau^*$. If $m(B)=\lambda_0(B)\leq m_0$ then
    \begin{equation*}
        m(B)\leq m(\overline{212111112122}11111222212211112^*121121111112211111\overline{221211111212}):=b_1.
    \end{equation*}
    Numerically $b_1\approx m_0-6.00958\cdot10^{-11}$. Moreover this value belongs to $L^\prime$.

    If $m(B)=\lambda_0(B)\geq m_1$ then
    \begin{equation*}
        m(B)\geq m(\overline{212221211111}222212211112^*121121111112211\overline{111112122212}):=a_2.
    \end{equation*}
    Numerically $a_2\approx m_0 + 8.762984\cdot 10^{-11}$. Moreover this value belongs to $L^\prime$.
\end{lemma}

\begin{proof}
    Suppose $m(B)=\lambda_0(B)\leq m_0$. By the proof of \Cref{lem3:intruder_branch_1}, we see that either $\lambda_0(B)<j_0-10^{-11}$ or $B=2111112\tau^*2111112$. Assuming the latter, just use the forbidden words 111112121, 111122121112, 21212, 222121112 to obtain the inequality. Now since $m(\overline{212111112122})=\lambda_0(212111112^*122)=3.32812\dots<j_0$ is an even semi-symmetric word, one can show that the above Markov value is actually in $L^\prime$ because of \Cref{lem:flahive}. 
    
    Now assume $m(B)=\lambda_0(B)\geq m_1$. Again by the proof of \Cref{lem3:intruder_branch_1} we see that either $\lambda_0(B)>m_0+10^{-10}$ or $B=2111112\tau^*2111111$. Assuming the latter, just use the forbidden words 111112121, 21212, 222121112. To see that the value of $a_2$ belongs to $L^\prime$ use the same argument as before.
\end{proof}

\begin{corollary}\label{cor3:second_characterization}
    If $B\in\{1,2\}^\Z$ satisfies $m(B)=\lambda_0(B)\in(b_1,a_2)$, then $B=\overline{w}w^*w1221211$. 
\end{corollary}
\begin{proof}
By \Cref{lem3:local_uniqueness} we have that either $B=122w^*111$ or $B=\tau^*$. In the second case, by \Cref{lem3:intruder_branch_2} one would have $\lambda_0(B)\leq b_1$ or $\lambda_0(B)\geq a_2$, which are impossible. Hence $B=122w^*111$ and by \Cref{cor3:first_characterization} we obtain  $B=\overline{w}w^*w12\dots$. Since 21212 is forbidden, is easy to see that $\lambda_0(\overline{w}w^*w1221211)<m_0+7.3\cdot 10^{-11}<a_2$ is the only possible continuation of $B$ with value less than $a_2$.
\end{proof}

\begin{lemma}\label{lem3:intruder_branch_3}
Suppose $B\in\{1,2\}^\Z$ satisfies $B=\tau^*$ and $m(B)=\lambda_0(B)\leq m_1$. Then 
\begin{equation*}
    m(B)\geq m(\overline{121111121222}212211112^*12112111111\overline{221211111212}):= a_1.
\end{equation*}
Numerically $a_1=3.33415633\ldots \approx j_0+6.26\cdot10^{-8}$. Moreover this value belongs to $L^\prime$. 
\end{lemma}
\begin{proof}
    Just use the forbidden words $2w2$, 111112121, $2w1111$, 21212 and \Cref{lem:flahive}.
\end{proof}

\begin{corollary}\label{cor3:gapsM}
We have $M\cap(j_0,m_0)\subset M\cap[a_1,b_1]$. In particular $(j_0,a_1)$ and $(b_1,m_0)$ are maximal gaps of $M$.
\end{corollary}

\begin{corollary}\label{cor3:gapL}
The interval $(b_1,a_2)$ is a maximal gap of $L$. 
\end{corollary}
\begin{proof}
    If $B\in\{1,2\}^\Z$ is a periodic orbit such that $m(B)=\lambda_0(B)\in(b_1,a_2)$, then by \Cref{cor3:second_characterization} one has $B=\overline{w}w^*w1221211$, which forces $B=\overline{w}$ and $m(B)=j_0<b_1$.
\end{proof}

Define $F$ to be the following set of words, including its transposes:
\begin{gather*}
    21212, 21112121, 2w2, 111112121, 111212111, 12111212, 2211112121122,  \\ 
    21221111212112112, 1212211112121121111, 221122111121211211212, \\
    2122111121211211112, w_r
\end{gather*}
where $w_r:=22w11$ is the self-replicating word. Define $F_0$ as the following set of words and their transposes:
\begin{gather*}
    1\tau, \tau1, 22\tau2, 12\tau22, 2112\tau211, 112\tau2112, 11112\tau211112, 211112\tau21111, \\
    1111112\tau211111, 2111112\tau211111
\end{gather*}

The proof of the following explicit characterization theorem is done in the same way as in the classical regions of \Cref{tab:knownregions}, but we include it for completeness.

\begin{theorem}\label{thm3:characterization1}
    We have that 
    \begin{equation*}
        M\cap(b_1,a_2)=(M\backslash L)\cap(b_1,a_2)=C\cup D\cup X
    \end{equation*}
    where
    \begin{multline*}
        C=\{\lambda_0\left(\overline{w}w^*w1221211\gamma\right):\gamma\in\{1,2\}^{\N}~\text{and $w1221211\gamma$} \\
        \text{does not contain any word from $F\cup F_0$}\}.
    \end{multline*}
    is a Cantor set and
    \begin{multline*}
        D=\Big\{\lambda_0\left(\overline{w}w^*w1221211\theta1121221\overline{w^T}\right):\theta\text{ finite word in 1,2 }, [0;\theta]\leq[0;\theta^T], \\
        \text{and $1221211\theta1121221$ does not contain any word from $F\cup F_0$}\Big\},
    \end{multline*}
    \begin{equation*}
        X=\left\{\lambda_0\left(\overline{w}w^*w122121121221\overline{w^T}\right),\lambda_0\left(\overline{w}w^*w1\overline{w^T}\right)\right\},
    \end{equation*}
    are sets of isolated points in $M$.
\end{theorem}

\begin{proof}
    Let $B\in\{1,2\}^\Z$ be such that $m(B)=\lambda_0(B)\in(b_1,a_2)$. By \Cref{cor3:second_characterization} we have that $B=\overline{w}w^*w1221211\gamma$ for some $\gamma\in\{1,2\}^\N$. Note that $\tilde{\gamma}=w1221211\gamma$ can not contain any word of $F\cup F_0$ except, possibly, the words $22w11$ or $11w^T11$. However if that tail $\tilde{\gamma}$ contains $22w11$, then \Cref{cor3:first_characterization} would imply that $\overline{w}$ contains $w12$. If $11w^T22$ also does not appear in that tail $\tilde{\gamma}$, then $\lambda_0(B)\in C$. By the same argument if $\lambda_0(B)\in C$, then again by \Cref{cor3:first_characterization}, \Cref{lem3:bounding_P_values} and the inequality $\lambda_0(\overline{w}w^*w1221211)<m_0+7.3\cdot 10^{-11}<a_2$ we also must have $m(B)=\lambda_0(B)\in M\cap(b_1,a_2)$.
    
    In case that $11w^T22$ does appear in $\tilde{\gamma}$, we will have that either \linebreak $B=\overline{w}w^*w1221211\theta1121221\overline{w^T}$ and so $\lambda_0(B)\in D$ or it is one of the finite possibilities of the set $X$ (some cases do not occur because they give Markov values outside $(b_1,a_2)$). As before one can show that $D\cup X\subset M\cap(b_1,a_2)$.

    Finally, to see that $D\cup X$ are isolated in $M$, given some $\lambda_0(B)\in D\cup X$, suppose that $m(B^{(n)})=\lambda_0(B^{(n)})\in(b_1,a_2)$ is any sequence converging to $\lambda_0(B)$. In particular one must have $B^{(n)}=\overline{w}w^*w1221211\dots$ for all large $n$. Since 
    $\lambda_0(B)-\lambda_0(B^{(n)})$ is the difference of just two continued fractions, the bi-infinite sequence $B^{(n)}$ must converge to $B$. In particular, for all large $n$, the subword $22w^T11$ will appear in $B^{(n)}$ which forces the rest of the continuation to be of the form $1\overline{w^T}$ (because of \Cref{cor3:first_characterization}). In particular one has $B^{(n)}=B$ and $m(B^{(n)})=\lambda_0(B)$ for all large $n$.
\end{proof}

\begin{corollary}
We have $p_{21211}^{(1)}=\max \left(M\cap(b_1,a_2)\right)$ and
\begin{equation*}
    p_{21211}^{(1)}=m(\overline{w}w^*w12212212112121111\overline{221211111212})\approx a_2-1.51\cdot10^{-11}.
\end{equation*}
\end{corollary}
\begin{proof}
    Just use \Cref{cor3:second_characterization} and the forbidden words 121211111, $22w11$ $2w2$, 21212, 12121112.
\end{proof}


\subsubsection{Separating $m_1$ from the intruder}

\begin{lemma}\label{lem3:intruder_branch_4}
    Suppose $B\in\{1,2\}^\Z$ satisfies, $B=\tau^*$, $m(B)<3.3341568$. Then $|m(B)-m_1|>3\cdot 10^{-11}$.
\end{lemma}

\begin{proof}
    Since $\lambda_0(1\tau^*)>m_1+10^{-7}$, $\lambda_0(\tau^*1)>m_1+10^{-7}$, $\lambda_0(12\tau^*21)<m_1-10^{-8}$ and $\lambda_0(22\tau^*22)>m_1+5\cdot 10^{-8}$ we have the two possibilities $12\tau^*22$ and $22\tau^*21$.
    \begin{itemize}
        \item Suppose we extend to $12\tau^*22$. Since $\lambda_0(112\tau^*22)>m_1+5\cdot 10^{-9}$ and $\lambda_0(212\tau^*222)<m_1-7\cdot10^{-9}$, we must continue as $212\tau^*221$. Since $\lambda_0(1212\tau^*221)<m_1-10^{-9}$ and $\lambda_0(212\tau^*2211)<m_1-5\cdot 10^{-10}$ we are led to $\lambda_0(2212\tau^*2212)>m_1+6\cdot 10^{-10}$.
        \item Suppose we extend to $22\tau^*21$. Since $\lambda_0(122\tau^*21)>m_1+5\cdot10^{-9}$ and $\lambda_0(22\tau^*211)>m_1+10^{-8}$ we must continue as $222\tau^*212$. Since $\lambda_0(222\tau^*2122)>m_1+7\cdot10^{-10}$ and $\lambda_0(1222\tau^*2121)<m_1-5\cdot 10^{-10}$ we must continue as $2222\tau^*2121$. Since 21212 is forbidden and $\lambda_0(12222\tau^*21211)>m_1+10^{-10}$ we must continue as $22222\tau^*21211$. Since $\lambda_0(22222\tau^*212112)>m_1+4\cdot 10^{-10}$ and $1111w^T22$ is forbidden we must continue as $\lambda_0(22222\tau^*2121111)>m_1+3\cdot10^{-11}$.
    \end{itemize}
\end{proof}

\begin{lemma}\label{lem3:intruder_branch_5}
    Suppose $B\in\{1,2\}^\Z$ satisfies $B=\tau^*$. If $m(B)=\lambda_0(B)\leq m_1$ then
    \begin{equation*}
        m(B)\leq m\left(\overline{211111212221}11212212212111221111212112111121211222\tau^*212112121111\overline{221211111212}\right) := b_2
    \end{equation*}
    If $m(B)=\lambda_0(B)\geq m_1$ then 
    \begin{equation*}
        m(B) \geq m\left(\overline{121111121222}2222212211112^*12112111111\overline{221211111212}\right) := a_3
    \end{equation*}
    and
    \begin{equation*}
        m(B)\leq m\left(\overline{111112122212}211112^*121121\overline{111112122212}\right)=3.334157259\ldots := b_3.
    \end{equation*}
    Numerically $b_2\approx m_1-5.27343\cdot10^{-10}$, $a_3\approx m_1+3.2499\dots10^{-11}$ and $b_3\approx j_1-2.22\cdot10^{-7}$. Moreover all this values belong to $L^\prime$.
\end{lemma}

\begin{proof}
    If $m(B)=\lambda_0(B)\leq m_1$, from the proof of \Cref{lem3:intruder_branch_4} we see that either $\lambda_0(B)<m_1-10^{-9}$ or $B=212\tau^*2211$ or $B=1222\tau^*2121$. If $B=212\tau^*2211$ then using the forbidden words from \Cref{lem3:local_uniqueness_forbidden_words} and \Cref{lem3:self_replication_forbidden_words} one has
    \begin{align*}
        m(B)\leq &m\left(\overline{211111212221}11212212\tau^*2211121221221211\overline{122212111112}\right) \\
        &\approx m_1-5.47542\cdot10^{-10},
    \end{align*}
    while if $B=1222\tau^*2121$ we have $m(B)\leq b_2$ given above. 

    If $m(B)=\lambda_0(B)\geq m_1$, then using the proof of \Cref{lem3:intruder_branch_4} and using that 21212, 111112121 and $2w2$ are forbidden, we have that $a_3\leq m(B)\geq b_3$ given above.
\end{proof}

\begin{corollary}\label{cor3:thirs_characterization}
    If $B\in\{1,2\}^\Z$ satisfies $m(B)=\lambda_0(B)\in(b_2,a_3)$, then $B=\overline{w}w^*w121112$. 
\end{corollary}
\begin{proof}
By \Cref{lem3:local_uniqueness} we have that either $B=122w^*111$ or $B=\tau^*$. In the second case, by \Cref{lem3:intruder_branch_5} one would have $\lambda_0(B)\leq b_2$ or $\lambda_0(B)\geq a_3$, which are impossible. Hence $B=122w^*111$ and by \Cref{cor3:first_characterization} we obtain  $B=\overline{w}w^*w12\dots$. Since 21212 is forbidden, is easy to see that $\lambda_0(\overline{w}w^*w121112)>m_1-7\cdot 10^{-11}>b_2$ is the only possible continuation of $B$ with values greater than $b_2$. All these values belong to $L^\prime$ because of \Cref{lem:flahive}.
\end{proof}

The proof of the following corollaries are analogous to that of \Cref{cor3:gapsM} and \Cref{cor3:gapL}

\begin{corollary}\label{cor3:gapsM2}
We have $M\cap(m_1,j_1)\subset M\cap[a_3,b_3]$. In particular $(m_1,a_3)$ and $(b_3,j_1)$ are maximal gaps of $M$.
\end{corollary}

\begin{corollary}\label{cor3:gapL2}
The interval $(b_2,a_3)$ is a maximal gap of $L$. 
\end{corollary}

Define $F_1$ as the set of following words including its transposes:
\begin{gather*}
    1\tau, \tau1, 22\tau22, 112\tau222, 2212\tau2212, 122\tau21, 22\tau211, 222\tau2122, 12222\tau21211, \\
    22222\tau212112, 22222\tau2121111
\end{gather*}

The proof of the following theorem is completely analogous to that of \Cref{thm3:characterization1}.

\begin{theorem}\label{thm3:characterization2}
    We have that 
    \begin{equation*}
        M\cap(b_2,a_3)=(M\backslash L)\cap(b_1,a_2)=C\cup D\cup X
    \end{equation*}
    where
    \begin{multline*}
        C=\{\lambda_0(\overline{w}w^*w121112\gamma):\gamma\in\{1,2\}^{\N}~\text{and $w121112\gamma$} \\
        \text{does not contain any word from $F\cup F_1$}\}.
    \end{multline*}
    is a Cantor set and
    \begin{multline*}
        D=\{\lambda_0(\overline{w}w^*w121112\theta211121\overline{w^T}):\theta\text{ finite word in 1,2 }, [0;\theta]\geq[0;\theta^T], \\
        \text{and $121112\theta211121$ does not contain any word from $F\cup F_1$}\},
    \end{multline*}
    \begin{equation*}
        X=\{\lambda_0(\overline{w}w^*w12111\overline{w^T})\},
    \end{equation*}
    are sets of isolated points in $M$.
\end{theorem}

\begin{corollary}
We have $p_{1112}^{(0)}=\min \left(M\cap(b_2,a_3)\right)$ and
\begin{equation*}
    p_{1112}^{(0)}=m\left(\overline{w}w^*w12111\overline{221211111212}\right)\approx b_2+4.54\cdot10^{-10}.
\end{equation*}
\end{corollary}


\section{Computations to bound the region with small dimension}

Let us denote by $t_{1/3}$ the positive real number such that  $t_{1/3}=\min\{t\in\R:D(t)=1/3\}$. In this short section we estimate this number:
\begin{theorem}
\begin{equation*}
    \nu_1=3.037002\dots<t_{1/3}<\mu_2=3.037311\dots,
\end{equation*}
where $\nu_1$ and $\mu_2$ are defined in \eqref{eq:nu_1} and \eqref{eq:mu_2}, respectively.
\end{theorem}

To prove this estimates, we will show that $D(\nu_1)<1/3<D(\mu_2)$.

There is maximal gap $(\mu_2,\nu_2)$ of the Markov spectra where 
\begin{equation}\label{eq:mu_2}
    \mu_2:=m\left(\overline{111221122211111112^*221122}\right)=3.037311\dots
\end{equation}
and
\begin{equation*}
    \nu_2:=m\left(\overline{11222111111112^*2211}\right)=3.03762983\dots
\end{equation*}

Indeed, suppose that $m(\underline{b})=\lambda_0(\underline{b})\in[\mu_2,\nu_2]$. Note that $\lambda_0(12^*1)>3.15$, $\lambda_0(2^*12)>3.06$, $\lambda_0(21112^*22)>3.04$, $\lambda_0(2111112^*221)>3.38$. so these words are forbidden. Since $\lambda_0(22^*2)<3$, we can assume $b_{-1}b_0^*b_1=12^*2$. Since 212 is forbidden and $\lambda_0(112^*21)<3.022$, we can assume $b_{-2}b_{-1}b_0^*b_1b_2=112^*22$. Since $\lambda_0(2112^*22)<3.01$,$\lambda_0(11112^*222)<3.035$ and 2111222, 212 are forbidden, we can assume that $b_{-4}b_{-3}b_{-2}b_{-1}b_0^*b_1b_2b_3b_4=11112^*2211$. Since $\lambda_0(211112^*2211)<3.034$ and 2111112221 is forbidden, we must continue as $b_{-6}b_{-5}b_{-4}b_{-3}b_{-2}b_{-1}b_0^*b_1b_2b_3b_4=1111112^*2211$.
\begin{itemize}
\item Assume that $b_5=1$. Since 2221112 is forbidden and $\lambda_0(21111112^*221111)<3.03725$, we must continue as $b_{-7}\dots b_6=11111112^*221111$. 

Since 
\begin{align*}
    \lambda_0(111222^*11111112)&>3.0379, \\
    \lambda_0(1111222^*111111111)&>3.03768, \\
    \lambda_0(221111222^*11111111221)&>3.03763, \\
    \lambda_0(221111222^*111111112222)&>3.0376299,
\end{align*}
these words are forbidden. Moreover, since
\begin{align*}
    \lambda_0(22221111222^*111111112221122)&>3.03763, \\
    \lambda_0(112221111222^*111111112221122)&>3.03762983,
\end{align*}
and 212, 121 are forbidden, we see that 222111122211111111222112 is also forbidden. Now we use the elementary fact that a continued fraction in $\{1,2\}$ that begins with $[2;2,2,1,1,1]$ and avoids 121, 212, 2221112, 1222111112, the previous four words and their transposes, is minimized with
\begin{equation*}
    [2;\overline{2,2,1,1,1,1,2,2,2,1,1,1,1,1,1,1,1,2}].
\end{equation*}
Similarly, since
\begin{align*}
    \lambda_0(11111222^*11111111)&>3.0377
\end{align*}
this word is forbidden. Finally, use that a continued fraction in $\{1,2\}$ that begins with $[2;1,1,1,1,1,1,1]$ and avoids all the previous forbidden words is minimized with
\begin{equation*}
    [2;\overline{1,1,1,1,1,1,1,1,2,2,2,1,1,1,1,2,2,2}].
\end{equation*}
From the previous inequalities one has that $m(\underline{b})=\lambda_0(\underline{b})\geq\nu_2$.

\item Assume $b_5=2$, which implies $b_6=2$ because 121 is forbidden. Using the fact that a continued fraction in $\{1,2\}$ that begins with $[0;2,2,1,1,2,2]$ and avoids 121, 212, 2111222, 2221112, 1222111112, 2111112221, 1111111222111, 1112221111111 is maximized with \[[0;\overline{2,2,1,1,2,2,1,1,1,2,2,1,1,2,2,2,1,1,1,1,1,1,1,2}]\] while a continued fraction in $\{1,2\}$ that begins with $[0;1,1,1,1,1,1]$ and avoids the same words is maximized with \[[0;\overline{1,1,1,1,1,1,1,2,2,2,1,1,2,2,1,1,1,2,2,1,1,2,2,2}],\] so we have that $\lambda_0(1111112^*221122)\leq\mu_2$.
\end{itemize}

From the previous computation we can conclude that
\begin{multline*}
    \Sigma_{\mu_2}=\Sigma_{\nu_2}=\Sigma\big(\{121, 212, 2111222, 2221112, 1222111112, 2111112221, \\
    1111111222111, 1112221111111\}\big).
\end{multline*}
Using the algorithms for computing the Hausdorff dimension \cite{MMPV} (the code is available at \url{https://github.com/Polevita/Gauss_Cantor_sets}), one verifies that $\dim_H(\Sigma_{\mu_2})=0.3367\dots>1/3$. Thus $t_{1/3}<\mu_2$.

Similarly, there is a maximal gap $(\mu_1,\nu_1)$ of Markov spectra where

\begin{align*}
    \mu_1 &:= m(\overline{22211221112211111111}112^*22112221122211\overline{11111111221112211222}) \\
    &=3.03699\dots
\end{align*}
and
\begin{equation}\label{eq:nu_1}
    \nu_1 := m(\overline{1112222211}222^*1111111111222\overline{1122222111})=3.037002\dots
\end{equation}
Moreover one has
\begin{multline*}
    \Sigma_{\mu_1}=\Sigma_{\nu_1}=\Sigma\big(\{121, 212, 2111222, 2111112221, 11122211111, 211111112221, \\ 
    12221111111112, 111111111112221, 12211222111111111, \\
    222211222111111111, \text{and their transposes}\}\big).
\end{multline*}
The program for computing the Hausdorff dimension gives that $\dim_H(\Sigma_{\mu_1})=0.3323\dots<1/3$. Thus $\nu_1<t_{1/3}$.


\appendix 

\section{More examples of $M\setminus L$ with intruder sets found numerically}\label{sec:further_examples}

Here, we present new non semisymmetric odd orbits that give rise to previously unknown regions of $M\setminus L$ containing also intruder sets. These orbits have been rigorously verified, although we do not provide detailed information about the structure of these regions. Instead, we only report the periodic orbit $\overline{w}$ and its self-replication height $j_1(w)$ (see \Cref{def:self-replication-height}).

\subsection{Examples with intruders in the alphabet $\{1,2\}$}\label{subsec:intruder_alphabet_12}

\subsubsection{Period 7}
\textbf{Before $t_1$}

\par\noindent\rule{\textwidth}{0.4pt}
\begin{equation*}
    m\left(\overline{1112^*122}\right)=\frac{2\sqrt{2026}}{27}=3.3341562\dots
\end{equation*}
\begin{align*}
    j_1(w)&=m\left(\overline{111212221211}ww^*ww1211\overline{122212111112}\right) \\
    &=3.3341574\dots \\
    &\approx j_0 + 1.204\cdot10^{-6}
\end{align*}
\par\noindent\rule{\textwidth}{0.4pt}

\subsubsection{Period 9}
\textbf{Before $t_1$}

\par\noindent\rule{\textwidth}{0.4pt}
\begin{equation*}
    m\left(\overline{22112^*1221}\right)=\frac{2\sqrt{65026}}{155}=3.2903478809\dots
\end{equation*}
\begin{align*}
    j_1(w)&=m\left(\overline{2111}ww^*ww211\overline{221211211212}\right)\\
    &=3.2903478821\dots \\
    &\approx j_0 + 1.14\cdot10^{-9}
\end{align*}
\par\noindent\rule{\textwidth}{0.4pt}

\textbf{After $t_1$}

\par\noindent\rule{\textwidth}{0.4pt}
\begin{equation*}
    m\left(\overline{11212^*1222}\right)=\frac{2\sqrt{63505}}{147}=3.4285984233\dots
\end{equation*}

\begin{align*}
    j_1(w)&=m\left(\overline{2122122121211121}122121211111212122112www^*w221212112212121111\overline{1222121211211212}\right) \\
    &=3.4285984248\dots \\
    &\approx j_0 + 1.45\cdot10^{-9}
\end{align*}
\par\noindent\rule{\textwidth}{0.4pt}

\subsubsection{Period 11}
\textbf{Before $t_1$}

\par\noindent\rule{\textwidth}{0.4pt}
\begin{equation*}
    m\left(\overline{222112^*12212}\right)=\frac{2\sqrt{2205226}}{903}=3.289037290587\dots
\end{equation*}
\begin{align*}
    j_1(w)&=m\left(\overline{1112}ww^*ww2\overline{2122122121122211}\right) \\
    &=3.289037290588\dots\\
    &\approx j_0 + 9.94\cdot10^{-13}
\end{align*}

\par\noindent\rule{\textwidth}{0.4pt}
\begin{equation*}
    m\left(\overline{122112^*12222}\right)=\frac{2\sqrt{2202257}}{901}=3.2941183949549\dots
\end{equation*}
\begin{align*}
    j_1(w)&=m\left(\overline{121221}2222121122ww^*ww12112212221211222212112\overline{2111}\right) \\
    &=3.2941183949559\dots \\
    &\approx j_0+1.01\cdot10^{-12}
\end{align*}

\par\noindent\rule{\textwidth}{0.4pt}
\begin{equation*}
    m\left(\overline{211112^*12222}\right)=\frac{\sqrt{1052677}}{309}=3.3203899266321\dots
\end{equation*}
\begin{align*}
    j_1(w)&=m\left(\overline{11212211}2121122ww^*ww\overline{2212111122}\right)\\
    &=3.3203899266366\dots \\
    &\approx j_0+4.5\cdot10^{-12}
\end{align*}

\par\noindent\rule{\textwidth}{0.4pt}
\begin{equation*}
    m\left(\overline{111112^*12221}\right)=\frac{\sqrt{848245}}{277}=3.32491758682\dots
\end{equation*}
\begin{align*}
    j_1(w)&=m\left(\overline{1111121222111112122112121121}ww^*ww\overline{1212112212111112221211111121}\right) \\
    &=3.32491758693\dots \\
    &\approx j_0+1.11\cdot10^{-10} 
\end{align*}
\par\noindent\rule{\textwidth}{0.4pt}

\textbf{After $t_1$}

\par\noindent\rule{\textwidth}{0.4pt}
\begin{equation*}
    m\left(\overline{211112^*12112}\right)=\frac{\sqrt{804613}}{269}=3.33458077937\dots,
\end{equation*}
\begin{align*}
    j_1(w)&=m\left(\overline{121111121222}ww^*ww\overline{221211111212}\right) \\
    &=3.33458077951\dots \\
    &\approx j_0 + 1.37\cdot10^{-10}
\end{align*}

\par\noindent\rule{\textwidth}{0.4pt}

\begin{equation*}
    m\left(\overline{221112^*12211}\right)=\frac{2\sqrt{451585}}{403}=3.33499128559\dots,
\end{equation*}
\begin{align*}
    j_1(w)&=m\left(\overline{211111212221}1121221221112122112121111ww^*ww2212111221\overline{122212111112}\right) \\
    &=3.33499128561\dots \\
    &\approx j_0 + 2.42\cdot10^{-11}
\end{align*}

\par\noindent\rule{\textwidth}{0.4pt}

\begin{equation*}
    m\left(\overline{121112^*12211}\right)=\frac{\sqrt{
    793885}}{267}=3.3370870586\dots
\end{equation*}
\begin{align*}
    j_1(w)&=m\left(\overline{121222121111}ww^*ww12122221121211112212111211221211\overline{122212111112}\right) \\
    &=3.3370870588\dots \\
    &\approx j_0 + \cdot10^{-10}
\end{align*}

\par\noindent\rule{\textwidth}{0.4pt}

\subsection{Examples with intruders in the alphabet $\{1,2,3\}$}\label{subsec:intruder_alphabet_123}

\subsubsection{Period 11}

\par\noindent\rule{\textwidth}{0.4pt}
\begin{equation*}
    m(\overline{w^*})=m\left(\overline{222123^*32121}\right)=\frac{\sqrt{65950645}}{2215}=3.6663657996727248\dots
\end{equation*}
\begin{align*}
    j_1(w)&=m\left(\overline{12}33212122121www^*w22211\overline{12}\right) \\
    &=3.6663657996727266 \\
    &\approx j_0 + 1.76\cdot10^{-15}
\end{align*}
\par\noindent\rule{\textwidth}{0.4pt}

\begin{equation*}
    m(\overline{w^*})=m\left(\overline{211233^*22222}\right)=\frac{2\sqrt{43917130}}{3573}=3.70948786758965611\dots
\end{equation*}
\begin{align*}
    j_1(w)&=m\left(\overline{21}12222ww^*ww2112212111122\overline{21}\right) \\
    &=3.70948786758965637\dots \\
    &\approx j_0 + 2.54\cdot10^{-16}
\end{align*}

\par\noindent\rule{\textwidth}{0.4pt}
\begin{equation*}
    m(\overline{w^*})=m\left(\overline{221233^*22122}\right)=\frac{2\sqrt{39803482}}{3393}=3.71883293796160443\dots
\end{equation*}
\begin{align*}
    j_1(w)&=m\left(\overline{21}12122ww^*ww2212212111\overline{12}\right) \\
    &=3.71883293796160474\dots \\
    &\approx j_0 + 3.11\cdot10^{-16}
\end{align*}

\par\noindent\rule{\textwidth}{0.4pt}
\begin{equation*}
    m(\overline{w^*})=m\left(\overline{111223^*22221}\right)=\frac{\sqrt{18671045}}{1127}=3.8340731702358434\dots
\end{equation*}
\begin{align*}
    j_1(w)&=m\left(\overline{21}1232321212112221www^*w11121\overline{12}\right) \\
    &=3.834073170235867\dots \\
    &\approx j_0 + 2.37\cdot10^{-14}
\end{align*}
\par\noindent\rule{\textwidth}{0.4pt}

\subsubsection{Period 13}

\par\noindent\rule{\textwidth}{0.4pt}
\begin{equation*}
    m(\overline{w^*})=m\left(\overline{1121233^*212212}\right)=\frac{\sqrt{432432029}}{5671}=3.666901797950443553\dots,
\end{equation*}
\begin{align*}
    j_1(w)&=m\left(\overline{21}112212ww^*ww112122121233\overline{21}\right) \\
    &=3.666901797950443594\dots \\
    &\approx j_0 + 4.09\cdot10^{-17}
\end{align*}

\par\noindent\rule{\textwidth}{0.4pt}
\begin{equation*}
    m(\overline{w^*})=m\left(\overline{1222123^*333212}\right)=\frac{\sqrt{1961781265}}{12062}=3.67202785701282708582\dots,
\end{equation*}
\begin{align*}
    j_1(w)&=m\left(\overline{21}12123332121233212www^*w122211\overline{12}\right) \\
    &=3.67202785701282708595\dots \\
    &\approx j_0 + 1.24\cdot10^{-19}
\end{align*}

\par\noindent\rule{\textwidth}{0.4pt}
\begin{equation*}
    m(\overline{w^*})=m\left(\overline{1111233^*211222}\right)=\frac{2\sqrt{58018690}}{4137}=3.68237856690483015\dots
\end{equation*}
\begin{align*}
    j_1(w)&=m\left(\overline{21}111222ww^*ww1111221\overline{212333}\right) \\
    &=3.68237856690483029\dots \\
    &\approx j_0 + 1.43\cdot10^{-16}
\end{align*}

\par\noindent\rule{\textwidth}{0.4pt}
\begin{equation*}
    m(\overline{w^*})=m\left(\overline{2123333^*211122}\right)=\frac{2\sqrt{978688657}}{16989}=3.6828536130416870816\dots
\end{equation*}
\begin{align*}
    j_1(w)&=m\left(\overline{21}111122ww^*ww2123321\overline{212333}\right) \\
    &=3.6828536130416870821\dots \\
    &\approx j_0 + 5.04\cdot10^{-19}
\end{align*}

\par\noindent\rule{\textwidth}{0.4pt}
\begin{equation*}
    m(\overline{w^*})=m\left(\overline{2121123^*321111}\right)=\frac{\sqrt{23677957}}{1321}=3.6835731285040127\dots
\end{equation*}
\begin{align*}
    j_1(w)&=m\left(\overline{12}211123332121221111www^*w212111\overline{12}\right) \\
    &=3.6835731285040135\dots \\
    &\approx j_0 + 8.62\cdot10^{-16}
\end{align*}

\par\noindent\rule{\textwidth}{0.4pt}
\begin{equation*}
    m(\overline{w^*})=m\left(\overline{1221123^*332121}\right)=\frac{2\sqrt{224820037}}{8127}=3.6899224888265465\dots
\end{equation*}
\begin{align*}
    j_1(w)&=m\left(\overline{21}1121ww^*ww1233\overline{21}\right) \\
    &=3.6899224888265477\dots \\
    &\approx j_0 + 1.18\cdot10^{-15}
\end{align*}

\par\noindent\rule{\textwidth}{0.4pt}
\begin{equation*}
    m(\overline{w^*})=m\left(\overline{2112223^*321211}\right)=\frac{\sqrt{474847685}}{5873}=3.71036950311272098\dots
\end{equation*}
\begin{align*}
    j_1(w)&=m\left(\overline{21}1111121221211www^*w211221\overline{12}\right) \\
    &=3.71036950311272101\dots \\
    &\approx j_0 + 3.48\cdot10^{-17}
\end{align*}

\par\noindent\rule{\textwidth}{0.4pt}
\begin{equation*}
    m(\overline{w^*})=m\left(\overline{2211223^*333211}\right)=\frac{2\sqrt{1008888170}}{17069}=3.7217177357480328218\dots
\end{equation*}
\begin{align*}
    j_1(w)&=m\left(\overline{21}1111121233211www^*w221121\overline{12}\right) \\
    &=3.7217177357480328223\dots \\
    &\approx j_0 + 4.85\cdot10^{-19}
\end{align*}

\par\noindent\rule{\textwidth}{0.4pt}
\begin{equation*}
    m(\overline{w^*})=m\left(\overline{3332123^*212222}\right)=\frac{13\sqrt{27914153}}{18364}=3.74014375992592658459\dots
\end{equation*}
\begin{align*}
    j_1(w)&=m\left(\overline{21}1212121112222www^*w333211\overline{12}\right) \\
    &=3.74014375992592658461\dots \\
    &\approx j_0 + 2.21\cdot10^{-20}
\end{align*}

\par\noindent\rule{\textwidth}{0.4pt}
\begin{equation*}
    m(\overline{w^*})=m\left(\overline{1111123^*212121}\right)=\frac{5\sqrt{879469}}{1251}=3.748201779800285\dots
\end{equation*}
\begin{align*}
    j_1(w)&=m\left(\overline{21}12121112121www^*w111111\overline{12}\right) \\
    &=3.748201779800302\dots \\
    &\approx j_0 + 1.64\cdot10^{-14}
\end{align*}
\par\noindent\rule{\textwidth}{0.4pt}

\subsubsection{Period 15}

\par\noindent\rule{\textwidth}{0.4pt}
\begin{equation*}
    m(\overline{w^*})=m\left(\overline{12111123^*3321111}\right)=\frac{2\sqrt{328987045}}{9841}=3.6862107565423097\dots
\end{equation*}
\begin{align*}
    j_1(w)&=m\left(\overline{21}11111ww^*ww\overline{12}\right) \\
    &=3.6862107565423102\dots \\
    &\approx j_0 + 5.55\cdot10^{-16}
\end{align*}

\par\noindent\rule{\textwidth}{0.4pt}
\begin{equation*}
    m(\overline{w^*})=m\left(\overline{11233223^*3332111}\right)=\frac{\sqrt{220820107229}}{126541}=3.71353948525976631184451\dots
\end{equation*}
\begin{align*}
    j_1(w)&=m\left(\overline{21}123332111212332111www^*w1123321\overline{12}\right) \\
    &=3.71353948525976631184467\dots \\
    &\approx j_0 + 1.61\cdot10^{-22}
\end{align*}

\par\noindent\rule{\textwidth}{0.4pt}
\begin{equation*}
    m(\overline{w^*})=m\left(\overline{11212223^*3321112}\right)=\frac{\sqrt{13368678133}}{31113}=3.716227943859402709\dots
\end{equation*}
\begin{align*}
    j_1(w)&=m\left(\overline{21}11112ww^*ww1121\overline{12}\right) \\
    &=3.716227943859402714\dots \\
    &\approx j_0 + 5.46\cdot10^{-18}
\end{align*}

\par\noindent\rule{\textwidth}{0.4pt}
\begin{equation*}
    m(\overline{w^*})=m\left(\overline{12121223^*3211121}\right)=\frac{\sqrt{2926052653}}{14547}=3.71849866205907477012\dots
\end{equation*}
\begin{align*}
    j_1(w)&=m\left(\overline{21}111212211121www^*w1212121\overline{12}\right) \\
    &=3.71849866205907477104\dots \\
    &\approx j_0 + 9.22\cdot10^{-19}
\end{align*}

\par\noindent\rule{\textwidth}{0.4pt}
\begin{equation*}
    m(\overline{w^*})=m\left(\overline{12212123^*2122222}\right)=\frac{\sqrt{17045652485}}{34949}=3.735700592729942052166\dots
\end{equation*}
\begin{align*}
    j_1(w)&=m\left(\overline{21}1122222ww^*ww12212111212123212223332\overline{21}\right) \\
    &=3.735700592729942052193\dots \\
    &\approx j_0 + 2.71\cdot10^{-20}
\end{align*}

\par\noindent\rule{\textwidth}{0.4pt}
\begin{equation*}
    m(\overline{w^*})=m\left(\overline{21222333^*2333211}\right)=\frac{\sqrt{1004912017213}}{268183}=3.73794386669548444907\dots
\end{equation*}

\begin{align*}
    j_1(w)&=m\left(\overline{21}111ww^*ww2111\overline{12}\right) \\
    &=3.73794386669548444912\dots \\
    &\approx j_0 + 5.68\cdot10^{-20}
\end{align*}

\par\noindent\rule{\textwidth}{0.4pt}
\begin{equation*}
    m(\overline{w^*})=m\left(\overline{22122123^*2122332}\right)=\frac{2\sqrt{35591840965}}{100911}=3.7390968279503780246669\dots
\end{equation*}
\begin{align*}
    j_1(w)&=m\left(\overline{21}1112121211122332www^*w2212211\overline{12}\right) \\
    &=3.7390968279503780246673\dots \\
    &\approx j_0 + 3.89\cdot10^{-22}
\end{align*}

\par\noindent\rule{\textwidth}{0.4pt}
\begin{equation*}
    m(\overline{w^*})=m\left(\overline{22212123^*2111212}\right)=\frac{\sqrt{3435483773}}{15647}=3.74595769375101444804\dots
\end{equation*}
\begin{align*}
    j_1(w)&=m\left(\overline{21}1111212ww^*ww222121112121\overline{12}\right) \\
    &=3.74595769375101444871\dots \\
    &\approx j_0 + 6.72\cdot10^{-19}
\end{align*}

\par\noindent\rule{\textwidth}{0.4pt}
\begin{equation*}
    m(\overline{w^*})=m\left(\overline{12212123^*2111221}\right)=\frac{2\sqrt{403326890}}{10721}=3.7464788778839061221\dots
\end{equation*}
\begin{align*}
    j_1(w)&=m\left(\overline{21}1111221ww^*ww122121112121\overline{12}\right) \\
    &=3.7464788778839061251\dots \\
    &\approx j_0 + 3.04\cdot10^{-18}
\end{align*}
\par\noindent\rule{\textwidth}{0.4pt}

\section{Special non semisymmetric orbits}\label{sec:special_orbits}

\subsection{Non semisymmetric orbits attaining the Markov value at two positions}

For the following orbits, we have marked their minimal periods at two positions $w_1\dots w_i^*\dots w_j^*\dots w_n$, $i\neq j$ to indicate the positions where the Markov value is attained.

\begin{equation*}
    m\left(\overline{211112^*22112211112^*2211221111122112}\right)=
    \frac{\sqrt{3440720805273180029}}{611699395}
\end{equation*}

\begin{equation*}
    m\left(\overline{11111112^*22222211112^*221122211222211111}\right)=\frac{10\sqrt{19523364336813387266}}{14571769931}
\end{equation*}

\begin{equation*}
    m\left(\overline{2211112^*2222112211112^*2222112211221122211}\right)=\frac{2\sqrt{72387925845470382258226}}{177722939265}
\end{equation*}

\begin{equation*}
    m\left(\overline{21112211221122111112^*2221122111112^*222112}\right)=\frac{\sqrt{38369319479152946329}}{2042162980}
\end{equation*}

\begin{equation*}
    m\left(\overline{1111112^*2211221111112^*2211221111111221122}\right)=\frac{\sqrt{276552940012055759845}}{5476905885}
\end{equation*}

\subsection{Non semisymmetric orbits with the same Markov value}

\subsubsection{First pair of orbits}

Define the following words of length 39
\begin{align*}
    u&:=112211222211221111122221122111221111122, \\ u^*&:=11221122221122111112^*2221122111221111122, \\ 
    v&:=112222112211221111122221122111112211122, \\
    v^*&:=11222211221122111112^*2221122111112211122.
\end{align*}

The periodic orbits $\overline{u}, \overline{u^T}$ and $\overline{v},\overline{v^T}$ are non semisymmetric, have minimal period 39 and are all \emph{distinct}. They have the remarkable property that they have the same Markov value:
\begin{equation*}
    j_0:=m\left(\overline{u}\right)=m\left(\overline{v}\right)=\frac{\sqrt{6486004003140082688401}}{26551371340}=3.0332042058\dots
\end{equation*}
Moreover, they both attain their Markov values at exactly one position in the minimal period, which is the position we have marked in $u^*$ and $v^*$.

These words are related by the following curious identities. Letting
\begin{equation*}
    w_l := 1122, \quad w_c := 2211221111122221122111, \quad w_r := 11122,
\end{equation*}
note that we have
\begin{equation*}
    u = w_l1122w_c2211w_r, \quad v = w_l2211w_c1122w_r.
\end{equation*}

\subsubsection{Second pair of orbits}\label{subsec:first_orbit_two_values}

Define the following words of length 39
\begin{align*}
    u&=122222112221122111122222112211112211221, \\
    u^*&=12222211222112211112^*2222112211112211221, \\
    v &= 122211222221122111122222112211221111221, \\
    v^* &= 12221122222112211112^*2222112211221111221.
\end{align*}
The periodic orbits $\overline{u}, \overline{u^T}$ and $\overline{v},\overline{v^T}$ are non semisymmetric, have minimal period 39 and are all \emph{distinct}. They have the remarkable property that they have the same Markov value:
\begin{equation*}
    j_0:=m\left(\overline{u}\right)=m\left(\overline{v}\right)=\frac{\sqrt{72400312134383876121601}}{88869071820}=3.027746947835738\dots
\end{equation*}
Moreover, they both attain their Markov values at exactly one position in the minimal period, which is the position we have marked in $u^*$ and $v^*$.

These words have a combinatorial structure somewhat similar to the example presented in \Cref{subsec:first_orbit_two_values}. Indeed, if we denote
\begin{equation*}
    w_l := 2222211222, \quad w_c := 112211112222211, \quad w_r := 221111221122,
\end{equation*}
then they have the following form
\begin{equation*}
    u = 1w_lw_cw_r1, \quad v = 1w_l^Tw_cw_r^T1.
\end{equation*}

\sloppy\printbibliography

\end{document}